\documentclass[12pt,twoside,a4paper]{article}

\usepackage{times}
\usepackage{amsmath,amssymb}
\usepackage{graphicx,color}

\newtheorem{theorem}{Theorem}[section]
\newtheorem{proposition}[theorem]{Proposition}
\newtheorem{lemma}[theorem]{Lemma}

\newlength{\espacioproof}
\newenvironment{proof}[1][Proof]{\noindent \textbf{#1.} }{\ \rule{.4em}{.4em} \vspace{\espacioproof}}

\newcommand{\R}{\mathbb{R}}
\newcommand{\N}{\mathbb{N}}

\renewcommand{\a}{\alpha}
\renewcommand{\b}{\beta}
\renewcommand{\O}{\Omega}
\newcommand{\e}{\varepsilon}
\newcommand{\g}{\gamma}
\renewcommand{\d}{\delta}
\newcommand{\D}{\Delta}
\newcommand{\s}{\sigma}

\renewcommand{\t}{\theta}
\renewcommand{\l}{\lambda}
\newcommand{\p}{\partial}
\newcommand{\f}{\varphi}

\newcommand{\inte}[1]{\overset{\circ}{#1}}

\DeclareMathOperator{\diam}{diam}
\DeclareMathOperator{\area}{area}
\DeclareMathOperator{\dist}{dist}
\DeclareMathOperator{\Card}{Card}
\DeclareMathOperator{\co}{co}
\DeclareMathOperator{\sgn}{sgn}

\newcommand{\ov}{\overline}
\newcommand{\un}{\underline}
\newcommand{\an}{\widehat}

\begin{document}
\title{Approximation of H\"older continuous homeomorphisms by piecewise affine homeomorphisms}
\author{Jos\'e C. Bellido and Carlos Mora-Corral}
\date{}
\maketitle

{\flushleft \small J.C.B.: ETSI Industriales, Universidad de Castilla-La Mancha. 13071 Ciudad Real. Spain. Email: JoseCarlos.Bellido@uclm.es \\
 C.M-C.: Mathematical Institute, University of Oxford. 24--29 St Giles'. Oxford OX1 3LB. United Kingdom. Email: mora-cor@maths.ox.ac.uk}

\subsection*{Abstract}

This paper is concerned with the problem of approximating a
homeomorphism by piecewise affine homeomorphisms. The main result is
as follows: every homeomorphism from a planar domain with a polygonal
boundary to $\R^2$ that is globally H\"older continuous of exponent $\a \in
(0,1]$, and whose inverse is also globally H\"older continuous of exponent
$\a$ can be approximated in the H\"older norm of exponent $\b$ by
piecewise affine homeomorphisms, for some $\b \in (0,\a)$ that only depends on $\a$. The proof
is constructive. We adapt the proof of simplicial approximation in
the supremum norm, and measure the side lengths and angles of the
triangulation over which the approximating homeomorphism is piecewise
affine. The approximation in the supremum norm, and a control on the
minimum angle and on the ratio between the maximum and minimum side
lengths of the triangulation suffice to obtain approximation in the
H\"older norm.

\section{Introduction}

This paper is concerned with the problem of approximating a homeomorphism by piecewise affine homeomorphisms.
As mentioned in Ball \cite{Ball01}, this problem arises naturally when one wants to approximate by finite elements the solution of a minimization problem in Nonlinear elasticity.
In that context, we are given a Lipschitz domain $\O \subset \R^n$ (typically, $n \in \{2,3\}$) and a function $h : \O \to \R^n$ that minimizes the elastic energy of a material, in a certain function space (typically, the Sobolev space $W^{1,p}$ for some $1 < p < \infty$).
In addition, in order for $h$ to be physically realistic, $h$ must be orientation-preserving and one-to-one (so as to avoid interpenetration of matter; see \cite{Ball81}).
Thus, every approximation of $h$ should also enjoy these two properties.
It is also pointed out in \cite{Ball01} that this question has theoretical interest too, since it would be a step towards generalizing Evans' \cite{Evans} result on the partial regularity of minimizers for integrands satisfying a certain growth condition.

When the original homeomorphism $h$ belongs to a Banach space $X$ of functions that includes piecewise affine functions, it is desirable to approximate $h$ by a piecewise affine homeomorphism both in the supremum norm and in the norm of $X$.
In fact, the difficulty of proving approximation of homeomorphisms $h$ by piecewise affine ones depends on the dimension $n$, the differentiability properties of $h$, and the norm in which this approximation is done.
In the context of Nonlinear elasticity explained above, ideally one assumes that $n \in \{2,3\}$, the fu‭nction $h$ and its inverse are in $W^{1,p}$, and looks for approximation in the $W^{1,p}$ norm.
Unfortunately, this is still an open problem, as put forward by Ball \cite{Ball01}.

The only positive results in this direction that we are aware of deal with approximation in the supremum norm.
In dimension $1$, the result that every homeomorphism can be approximated by a piecewise affine homeomorphism in the supremum norm is trivial.
The first proof in dimension $2$ seems to be Rad\'o's \cite{Rado} (see also Moise \cite{Moise} and Brown \cite{Brown}).
The result in dimension $3$ is due to Moise \cite{Moise52} and Bing \cite{Bing54}.
For dimensions $5$ and higher, the result for contractible spaces follows from theorems of Connell \cite{Connell}, Bing \cite{Bing63}, Kirby \cite{Kirby} and Kirby, Siebenmann and Wall \cite{KSW} (for a proof see, e.g., Rushing \cite{Rushing} or Luukkainen \cite{Luukkainen}).
Finally, Donaldson and Sullivan \cite{DS} proved that the result is false in dimension $4$.

In this paper we consider the problem of approximation in the H\"older norm and in dimension $2$.
Our main result reads as follows.
\begin{theorem}\label{th:introduction}
Let $\O \subset \R^2$ be a closed polygon. Let $0 < \a \leq 1$. Let
$h \in C^{\a} (\O , \R^2)$ be a homeomorphism such that $h^{-1} \in C^{\a} (h(\O) , \R^2)$.
Then there exists $0 < \b < \a$, depending only on $\a$, such that
for each $\e >0$ there exists a piecewise affine homeomorphism $f:
\O \to \R^2$ with $\| f - h \|_{\b} < \e$.
\end{theorem}
Here $C^{\a}$ denotes the Banach space of globally H\"older continuous functions of exponent $\a$, with norm $\|\cdot\|_{\a}$.

We follow the proof of Moise \cite{Moise}, where approximation in the supremum norm is proved, but there the construction is not explicit. In this paper we make an explicit construction, measure the lengths and angles of the triangulation, and show that if we have approximation in the supremum norm and a control on the angles of the triangulation and on the ratio between the maximum and minimum side lengths of the triangulation, then we have approximation in the H\"older norm. Interestingly, when $\a =1$, the constructed triangulation is regular in the sense of Ciarlet \cite{Ciarlet}.

We now describe the outline of this paper.
Section \ref{se:notation} introduces the notations and definitions that will be used throughout the paper.
Section \ref{se:plan} describes the plan of the proof of Theorem \ref{th:introduction},
and each of the remaining sections (\ref{se:skeleton}, \ref{se:extension}, \ref{se:triangulation} and \ref{se:estimates}) is devoted to a specific step of the proof.
Following the notation of Theorem \ref{th:introduction}, in Section \ref{se:skeleton} we show how to refine a given triangulation of $\O$, and how to construct a piecewise affine function over the skeleton of the refined triangulation that approximates $h$ in the supremum norm; we also measure the minimum and maximum lengths of the triangulation.
In Section \ref{se:extension} we extend an arbitrary homeomorphism defined on the boundary of a triangle to a homeomorphism defined on the whole triangle; we triangulate the original triangle and measure the lengths and the angles of the triangulation.
Section \ref{se:triangulation} constructs the piecewise affine homeomorphism $f$. The idea is as follows: we start with a fine regular triangulation of $\O$; then we add vertices in the skeleton and construct an approximating piecewise affine homeomorphism $g$ on the skeleton, using the result of Section \ref{se:skeleton}; then we extend this piecewise affine homeomorphism $g$ on the skeleton
to an approximating piecewise affine homeomorphism $f$ on the whole triangulation, using the result of Section \ref{se:extension}. The outcome of Section \ref{se:triangulation} is a piecewise affine homeomorphism $f$ that approximates $h$ in the supremum norm, and we also estimate the lengths and angles of the triangulation over which $f$ is piecewise affine.
In Section \ref{se:estimates} we show general a priori bounds in the H\"older norm of any piecewise affine function $u$, in terms of the lengths and angles of the triangulation over which $u$ is piecewise affine. Finally, we show how these a priori bounds demonstrate that the piecewise homeomorphism $f$ constructed in Section \ref{se:triangulation}  approximates $h$ also in the H\"older norm, thus concluding Theorem \ref{th:introduction}.

\section{Notations and definitions}\label{se:notation}

Two key concepts are used in the construction of this paper: a \emph{complex} and a \emph{piecewise affine function}.
We assume no previous knowledge of complexes; rather, they serve only as a useful notation.
For us, a \emph{complex} means what is usually referred to as a \emph{Euclidean finite complex in $\R^2$} (see, for example, Chapter 0 of Moise \cite{Moise}).
Specifically, a $1$-dimensional complex is a non-empty set $K$ such that:
\begin{itemize}
\item every element of $K$ is either a closed segment in $\R^2$ or a singleton;
\item $\{p\}\in K$ if and only if $p$ is an endpoint of a segment of $K$;
\item if $e,c$ are segments in $K$ and $e \cap c \neq \varnothing$, then $e \cap c$ has exactly one element, which is an endpoint of both $e$ and $c$.
\end{itemize}
A $2$-dimensional complex is a non-empty set $K$ such that:
\begin{itemize}
\item every element of $K$ is either a triangle in $\R^2$, or a closed segment or a singleton;
\item $\{p\}\in K$ if and only if $p$ is an endpoint of a segment of $K$;
\item a segment belongs to $K$ if and only if it is a side of a triangle of $K$;
\item if $\s,\tau$ are triangles in $K$ and $\s \cap \tau \neq \varnothing$, then $\s \cap \tau$ is either a singleton whose only member is a vertex of both $\s$ and $\tau$, or a segment which is a side of both $\s$ and $\tau$.
\end{itemize}
If $K$ is an $n$-dimensional complex (for some $n \in \{1,2\}$) then
the set $K^0$ is the set formed by the singletons of $K$,
the set $K^1$ is the set formed by the segments of $K$, and (if $n=2$)
$K^2$ is the set formed by the triangles of $K$.
Naturally, $K = K^0 \cup K^1$ for any $1$-dimensional complex $K$, and
$K = K^0 \cup K^1 \cup K^2$ for any $2$-dimensional complex $K$.
Of course, if $K$ is a $2$-dimensional complex then $K^0 \cup K^1$ is a $1$-dimensional complex, sometimes called the \emph{skeleton} of $K$.

Given a complex $K$, it is convenient to work with the set $\tilde{K}^0$ defined as follows: if $K^0 = \{ \{p_1\}, \ldots , \{p_n\} \}$ for some $n \in \N$ and some $p_1, \ldots, p_n \in \R^2$ then
$\tilde{K}^0 := \{ p_1, \ldots , p_n \}$. And reciprocally, $\tilde{K}^0$ defines $K^0$ univocally. Of course, $\bigcup K^0 = \tilde{K}^0$.

If $K$ is a $2$-dimensional complex we call $K$ a \emph{triangulation} of $\bigcup K$.

In the previous paragraphs we have mentioned segments. Although the definition should be known, we make it precise.
Given two different points $x, y$ in $\R^2$, we will denote $[x, y]$ the segment in $\R^2$ with endpoints $x$ and $y$ equipped with the total order given by the bijection
\[
 \begin{array}{rcl}
 [0,1] & \to & [x, y] \\
 t & \mapsto & x + t (y - x) .
 \end{array}
\]
The set $(x,y)$ equals $[x,y] \setminus \{x,y\}$ and inherits the total order.
The set $\ov{xy}$ equals $[x,y]$ as a set, but with no order structure.

Another key concept of this paper is the one of \emph{piecewise affine function}.
Let $K$ be a complex. A function $f: \bigcup K \to \R^2$ is \emph{piecewise affine over $K$} when $f|_{\s}$ is affine for all $\s \in K$. Since every element of $K$ is closed, this $f$ is automatically continuous.
Now let $I$ be a totally ordered set, and choose $n\in \N$ points
\[
 a_1 < \cdots < a_n
\]
of $I$. A function $f: [a_1,a_n] \to \R^2$ is \emph{piecewise affine over $\{ a_1 , \ldots , a_n \}$} when $f|_{[a_i, a_{i+1}]}$ is affine for all $1\leq i\leq n-1$. Since $[a_i, a_{i+1}]$ is closed for all $1\leq i\leq n-1$, this $f$ is automatically continuous.

Another elementary concept, which nevertheless deserves some care in the notation, is that of \emph{polygon}.
Let $n\geq 3$ be a natural number. Consider $n$ points $a_1 , \ldots, a_n \in \R^2$.
We will say that $a_1 \cdots a_n$ is a \emph{well-defined polygon} if there exists a homeomorphism $f$ from the unit circle of $\R^2$ onto the set
\begin{equation}\label{eq:Jordancurve}
 \ov{a_n a_1} \cup \bigcup_{k=1}^{n-1} \ov{a_k a_{k+1}}
\end{equation}
and there exist $0 \leq t_1 < \cdots < t_n < 2 \pi$ such that $f(\cos t_i , \sin t_i) = a_i$ for all $i \in \{1, \ldots, n\}$.
In this case, $a_1 \cdots a_n$ equals the set \eqref{eq:Jordancurve} union the connected component of the complement of \eqref{eq:Jordancurve}.
We will say that $a_1 \cdots a_n$ is a \emph{well-defined $n$-gon} when it is a well-defined polygon and there do not exist $b_1 , \ldots, b_{n-1} \in \R^2$ such that
\[
 \ov{a_n a_1} \cup \bigcup_{k=1}^{n-1} \ov{a_k a_{k+1}} = \ov{b_{n-1} b_1} \cup \bigcup_{k=1}^{n-2} \ov{b_k b_{k+1}} .
\]
Of course, a $3$-gon is called a triangle, and a $4$-gon, a quadrilateral.
A \emph{closed polygon} is a compact set in $\R^2$ that coincides with the closure of its interior, and whose boundary is a finite union of Jordan curves, each of them is the boundary of a well-defined polygon.

If $A\subset \R^2$ then $\inte{A}$ denotes the interior of $A$ in the topology of $\R^2$, except if $A$ is a segment, in which case $\inte{A}$ denotes the set $A$ minus its endpoints.

Balls follow the usual notation: $\bar{B}(a,r)$ is the closed ball centred at $a\in \R^2$ with radius $r>0$;
and $B(A,r)$ and $\bar{B}(A,r)$ are the open and closed neighbourhoods, respectively, of $A\subset
\R^2$ with radius $r>0$. We will always use the Euclidean distance.

Let $\O \subset \R^2$ be compact.
Let $\|\cdot\|_{\infty}$ denote the supremum norm on $\O$. The supremum norm on a subset $S$ of $\O$ will be denoted by $\|\cdot\|_{\infty,S}$. For each $0 < \a \leq 1$, let $|\cdot|_{\a}$ denote the H\"older seminorm, and $\|\cdot\|_{\a}$ the H\"older norm, both of exponent $\a$. The Banach space of globally H\"older continuous functions of exponent $\a$ from $\O$ to $\R^2$ is denoted by $C^{\a}( \O , \R^2)$.

\section{Plan of the proof}\label{se:plan}

In this section we describe the main lines of our construction of a piecewise affine homeomorphism that approximates a given homeomorphism, in the conditions of Theorem \ref{th:introduction}.
We believe that this description may serve as a guide in order to enjoy an easier reading of the rest of
the paper.

We follow the construction due to Moise \cite{Moise} in the two-dimensional case.
It is proved there that given a closed polygon $\O$, a homeomorphism $h\in C(\O,\R^2)$, and $\e >0$, there exist a
triangulation $K$ of $\O$, and a piecewise affine homeomorphism $f : \O \to \R^2$ such that $\| h-f\|_{\infty} \leq \e$.
Extending Moise's construction to other functional spaces (for instance, Sobolev spaces $W^{1,p}$ or H\"older spaces $C^{\a}$) is a very delicate issue, since in his construction there is no control at all on the gradient of the approximation.
Indeed, if we think in terms of the mathematical theory of finite element approximation by piecewise affine functions in the Sobolev norm (see, e.g., Ciarlet \cite{Ciarlet}), we find two major difficulties:
\begin{itemize}
\item In order to guarantee the injectivity of $f$, we cannot use nodal values; that is to say, in general,
$f(a)\neq h(a)$ for every vertex $a$ of the triangulation $K$. This occurs in Moise's construction too.
The problem here is that if we use nodal values then the orientation of the triangles may change and, hence, $f$ will, in general, fail to be one-to-one. Of course, this problem would dissapear if the original function $h$ were a diffeomorphism.

\item Except when $\a = 1$, the triangulation that we construct is not regular in the sense of Ciarlet (see, e.g., Ciarlet \cite{Ciarlet} or Zl\'amal \cite{Zlamal}). In other words, for every $\e>0$ we construct a triangulation $K_{\e}$ of $\O$, and, in our construction, the minimum angle of the triangles of the triangulation $K_{\e}$ tends to zero as $\e$ tends to zero. It is very well-known that the regularity of a triangulation (i.e., a positive lower bound independent of $\e$ for the minimum angle) is crucial in order to have approximation results.
\end{itemize}

Roughly speaking, what happens is that both the triangulation $K$ and the
approximation $f$ have to adapt themselves to the target function
$h$, if we want $f$ to be a homeomorphism. To circumvent these
difficulties, in order to prove Theorem \ref{th:introduction}, we will follow 
Moise's construction, but making everything completely explicit, so
that we can exactly estimate the norm $\|f-h\|_\b$ at the end. We
have to use fine analytical and geometrical arguments to prove Theorem \ref{th:introduction},
whereas Moise just needed topological arguments to get his results; in fact, his construction is not explicit.

Let $\a \in (0,1]$ be the exponent of H\"older continuity of $h$, and
$\tilde{\a} \in (0,1]$ the exponent of H\"older continuity of $h^{-1}$.
In the statement of Theorem \ref{th:introduction} we assumed $\a = \tilde{\a}$; however, allowing $\a$ and $\tilde{\a}$ to be different gains insight in the proof and provides slightly better estimates.

We fix an initial triangulation $M$ of $\O$. The proof of Theorem \ref{th:introduction} is structured in the following four steps:

\begin{description}
\item[{\it Step 1: Approximation in the skeleton.}]
Consider the skeleton $M^1$ of the triangulation $M$. In this step we
find a subdivision $K^1$ of $M^1$ and a piecewise affine function $f$ over the
subdivision such that $f$ is a homeomorphism on $\bigcup K^1 = \bigcup M^1$ and
\[
 \|f-h\|_{\bigcup M^1,\infty}\leq \e .
\] 
Section \ref{se:skeleton} is devoted to the proof this step. It is important to remark
that the function $f$ is not the piecewise affine interpolant of
$h$ (i.e., in general $f(a) \neq h(a)$ for $a\in K^0$), although $f(a)=h(a)$ for all $a\in M^0$, i.e., for any vertex $a$ of the original triangulation $M$.

The construction of the subdivision is explicit, in the sense that we estimate the length of the edges in $K^1$ in terms of
$\e$. Precisely, for each $e\in K^1$ the inequalities 
\[
 B_1 \e^{b_1} \leq |e| \leq B_2 \e^{b_2}
\]
hold, where the exponents $b_1 \geq b_2 \geq 1$ depend only on $\a, \tilde{\a}$, and the constants $B_1,B_2 >0$ depend only on $M$, $h$ and $h^{-1}$. Constants and exponents are calculated explicitly.

\item[{\it Step 2: Extension of a homeomorphism from the boundary of a triangle to the whole triangle.}]
For every triangle $\D \subset \R^2$ and every piecewise affine homeomorphism
$f : \p\D \to \R^2$, there exists a piecewise affine homeomorphism $\tilde{f} : \D \to \R^2$ that extends $f$.
The existence of such an extension is known as the \emph{piecewise affine} (or \emph{PL}) \emph{Schoenflies
Theorem}. The proof in Moise \cite{Moise} is constructive, but of course there are no estimates on the triangulation
parameters, since they are not needed in that topological context. Instead, we give a different constructive proof of the result,
based on an idea borrowed from Gupta and Wenger \cite{GW}, according to which we perturb the triangle $\D$ to get a well-defined convex $w$-gon, where $w$ is the number of vertices of the triangulation of $\p \D$ over which $f$ is piecewise affine. We estimate all the parameters (side lengths and angles) of the constructed triangulation of $\D$. Section \ref{se:extension} is devoted to this step, and Theorem \ref{th:Gupta} is the main result in this section.

\item[{\it Step 3: Construction of the triangulation.}]
We put steps 1 and 2 together in order to build the
final triangulation and the piecewise affine homeomorphism. We proceed in the following way: we start with a regular and
sufficiently fine triangulation $M$. By step 1, we find a
subdivision $K^1$ of $M^1$ and a piecewise affine homeomorphism $f$ on
the subdivision. Now we use step 2 to extend the
homeomorphism $f$ (initially defined only on the skeleton of the triangulation) to the whole $\O$; we do this by finding a triangulation of each triangle of $M$.
In this way, we get a new triangulation $K$ of $\O$ and an extension of $f$ (denoted again by $f$)
such that:
\begin{enumerate}
\item $K$ is the union of the triangulations of the triangles of $M$;
\item $f$ is piecewise affine over the triangulation $K$;
\item $f$ is a homeomorphism (since it is a homeomorphism over each triangle of $M$ and on the skeleton of th triangulation).
\end{enumerate}

This step is standard following Moise's proof, again except for the estimates of the triangulation parameters. Thus, we obtain a triangulation $K$ of $\O$ such that
\begin{itemize}
\item $\sin \f \geq A_0 \e^{a_0}$, for all angles $\f$ of all triangles of the triangulation $K$;
\item $A_1 \e^{a_1} \leq |e|\leq A_2 \e^{a_2}$, for all sides $e\in K^1$,
\end{itemize}
and $\|h-f\|_\infty \leq \e$.
The quantities $A_0, A_1, A_2, a_0, a_1, a_2$ are calculated explicitly. 
The exponents $a_0, a_1, a_2$ depend only on $\a, \tilde{\a}$, and they satisfy
$a_0 \geq 0$ and $a_1 \geq a_2 \geq 1$. The constants $A_0, A_1, A_2 > 0$ depend also on $h$, $h^{-1}$ and $\O$.
Section \ref{se:triangulation} is devoted to this step, and Theorem \ref{th:constructiontriangulation} is the
main result in that section.

We would like to remark that everything that we have done up to now has
followed essentially the ideas of the classical proof, but with
the important difference that we estimate all the parameters of the
triangulation in terms of powers of $\e$.
The arguments of the following step are new, and allow us to
extend the classical approximation result to H\"older spaces.

\item[{\it Step 4: From $L^{\infty}$ estimates to $C^{\b}$ estimates.}]
Section \ref{se:estimates} is devoted to the final step of obtaining approximation in the H\"older norm.
Up to now, our construction gives a triangulation
$K$ and a piecewise affine homeomorphism $f$ that approximates
$h$ in the supremum norm. In this step now we prove that actually $f$ approximates $h$
in the H\"older norm $C^\b$; this $\b$ depends only on $\a, \tilde{\a}$. This is a consequence of the
approximation on the supremum norm and a control on the
triangulation in terms of $\e$. We start using the
interpolation inequality
\[ 
 |u|_{\b} \leq 2 \|u\|_{\infty}^{1-\frac{\b}{\a}} |u|_{\a}^\frac{\b}{\a},
\]
valid for any $u \in C^{\a}$ and $0 < \b \leq \a \leq 1$.
Using the previous step we obtain that 
\[
 |h-f|_{\b} \leq 2 \e^{1-\frac{\b}{\a}} |h-f|_{\a}^{\frac{\b}{\a}},
\]
so that we just need a priori bounds on the seminorm $|h-f|_{\a}$.
Since in the construction of $f$ (steps 1, 2 and 3) no care has been taken in the gradients, the fact that $f$ is close to $h$ in the $C^{\a}$ norm is not guaranteed at all.
Hence the only valid estimate of $|h-f|_{\a}$ is the trivial one: $|h-f|_{\a} \leq |h|_{\a} + |f|_{\a}$.
Thus, we just need a priori bounds on $|f|_{\a}$.
Here we see why in our construction we cannot get approximation in the $C^{\a}$ norm.
This impossibility is not just because of the proof, but it is a general fact (see the comments at the end of the paper).

The a priori bounds on $|f|_{\a}$ are given in Propositions \ref{prop:piecewise} and \ref{prop:errorinterpolant}.
In Proposition \ref{prop:piecewise} we estimate the H\"older seminorm of any piecewise affine function $u$ over the triangulation $K$ in terms of the parameters of the triangulation (and, ultimately, of $\e$) and the $L^{\infty}$ norm of $u$.
In Proposition \ref{prop:errorinterpolant} we estimate the $C^{\a}$ seminorm of the interpolant over $K$ of any $C^{\a}$ function $u$, in terms of $C^{\a}$ norm of $u$, and of the parameters of the triangulation (so, ultimately, in terms of $\e$).
Now, using those estimates we get $|h-f|_{\a} \leq C \e^{a_4}$,
for some constant $C>0$ depending on $h$, $h^{-1}$ and $\O$, and some exponent $a_4 \leq 0$ depending only on $\a, \tilde{\a}$. As usual, both $C$ and $a_4$ are explicit.
Consequently, for every $0 < \b < \frac{\a}{1-a_4}$ we have
\[
 |h-f|_{\b} \leq D \e^{1-\frac{\b}{\a}(1-a_4)} .
\]
For some constant $D>0$ depending on $h$, $h^{-1}$ and $\O$. 
The fact that $1-\frac{\b}{\a}(1-a_4)$ is positive concludes the result.

\end{description}

We finish this section with some remarks about possible generalizations of Theorem \ref{th:introduction} to other function spaces.

Let $X$ be a Banach space continuously included in $C(\O,\R^2)$ and such that $X$ contains all piecewise affine functions from $\O$ to $\R^2$.
For example, $X$ can be a Sobolev space $W^{1,p}$ with $p>2$, or a space of uniformly continuous functions with prescribed modulus of continuity.
We believe that steps 1, 2 and 3 above can be adapted with only minor modifications, in the context of the function space $X$.
Step 4 can also be adapted easily, but in this case we would obtain approximation in the function space $Y$, where $Y$ is an interpolation space between $C(\O,\R^2)$ (or $L^{\infty}(\O,\R^2)$) and $X$.
This is the only reason why our construction does not work for proving $W^{1,p}$ approximation:
because for all $1 \leq p < q \leq \infty$, the space $W^{1,p}$ is not an interpolation space between $L^{\infty}$ and $W^{1,q}$.
In contrast, we believe that the construction of this paper can be easily adapted to get approximation in the fractional Sobolev space $W^{s,p}$, for some $0 < s < 1$ and $2 < p \leq \infty$.

\section{Approximation in the skeleton}\label{se:skeleton}

In this section we approximate a given homeomorphism on a $1$-dimensional complex $M$
by a homeomorphism piecewise affine over a refinement of $M$.

Since the proof of Theorem \ref{th:skeleton} below is long, we have decided first to explain its main ideas, without proofs.
This is done in the next paragraphs.
 
Let $M$ be a $1$-dimensional complex in $\R^2$.
Let $h : \bigcup M \to \R^2$ be a homeomorphism such that $h$ and $h^{-1}$ are H\"older continuous.
In Section \ref{se:triangulation}, $M$ will be the skeleton of a fine quasiuniform triangulation of $\O$, and $h$ will be the restriction of the homeomorphism from $\O$ to $\R^2$ that we want to approximate, but in this section we are only concerned with what happens in the skeleton.
We start by subdividing (refining) the $1$-dimensional complex $M$ in a uniform way to obtain a new $1$-dimensional complex $L$ such that $\diam h(e) \leq \e /3$ for all $e \in L^1$.
Thus, if we define $N_v := \bar{B} (h(v), \b)$ for each $v\in \tilde{L}^0$ and a suitable $\b \leq \e/3$, then we show that
\[
 N_u \cap N_v=\varnothing , \qquad u, v \in \tilde{L}^0 ,
\]
and moreover, if $v\in \tilde{L}^0$ does not belong to $e\in L^1$ then
\[
 N_v \cap \bar{B} (h(e),\frac{\e}{3}) = \varnothing .
\]
Fix $e\in L^1$, and choose one of the two possible orientations for $e$, so that $e = [\inf e, \sup e]$.
Accordingly, equip $h(e)$ with the natural order given by $e$ and the bijection $h$.
Let $x_e$ be the last point of the arc $h(e)$ that lies on $N_{\inf e}$, and $y_e$ the first point that follows $x_e$ and lies on $N_{\sup e}$.
We proceed by finding a fine enough uniform partition $\{ w_{e,0}, \ldots, w_{e,N_e} \}$ of the segment $[h^{-1}(x_e), h^{-1}(y_e)]$, and define $g_e : [h^{-1}(x_e), h^{-1}(y_e)] \to \R^2$ as the piecewise affine function over $\{ w_{e,0}, \ldots, w_{e,N_e} \}$ that coincides with $h$ in $\{ w_{e,0}, \ldots, w_{e,N_e} \}$.
It may well happen that the intersection of $g_e [h^{-1}(x_e),h^{-1}(y_e)]$ with $N_{\inf e}$ or $N_{\sup e}$ is not empty.
To avoid this possibility, which could yield the non-injectivity of the piecewise affine function that we are constructing, we do the following: define $p_e$ as the last point of $g_e [h^{-1}(x_e),h^{-1}(y_e)]$ that lies on $N_{\inf e}$, and
$q_e$ as the first point following $p_e$ that lies on $N_{\sup e}$.
Now from the partition $\{ w_{e,0}, \ldots, w_{e,N_e} \}$ we get the partition $\{u_{e,0}, \ldots, u_{e,m_e}\}$ of $[p_e,q_e]$ (for some $m_e \leq N_e$) defined by
\[
 u_{e,0} = p_e, \qquad u_{e,m_e} = q_e, \qquad \{ u_{e,0}, \ldots, u_{e,m_e} \} = [p_e,q_e] \cap \{ w_{e,0}, \ldots, w_{e,N_e} \} .
\]

As $h$ is a homeomorphism, then the points $g_e(u_{e,0}), \ldots, g_e(u_{e,m_e})$ are all different.
It is important to observe that the function $g|_{[p_e,q_e]}$ need not be injective, since a loop may be formed because of the geometry of the arc $h[p_e,q_e]$.
If this happens, we just remove the loop (see Figure \ref{fig:loop}), and it is easy to show (Lemma \ref{le:fginjective}) that there exists a new injective piecewise affine function
$f_e : [p_e,q_e] \to \R^2$ over a partition of at most $m_e + 1$ elements such that
\[
 f_e(p_e) = g_e(p_e), \qquad f_e(q_e) = g_e(q_e), \qquad f[p_e,q_e] \subset g_e[p_e,q_e] .
\]
It is important to remark that $g_e$ was obtained by using nodal values of the function $h$, but this is no longer the case for $f_e$; in other words, to ensure injectivity we have to take values in the vertices that are not necessarily nodal values of $h$.

Finally, we define $f: \bigcup M \to \R^2$ as the only piecewise affine function such that, for each $e\in L^1$,
\begin{align*}
 & f|_{[p_e,q_e]} = f_e , \\
 & f|_{[\inf e,p_e]} \quad \mbox{is affine with} \quad f(\inf e) = h(\inf e), \quad f(p_e) = g_e(p_e) , \\
 & f|_{[q_e,\sup e]} \quad \mbox{is affine with} \quad f(q_e) = g_e(q_e), \quad f(\sup e) = h(\sup e) .
\end{align*}
Then, $f$ is a homeomorphism, since it is injective on each $[p_e,q_e]$ (for $e \in L^1$), on each $[q_e,\sup e] \cup [u_c,p_c]$ (for all $e,c \in L^1$ such that $\sup e=\inf e$), and the interior of the images of these sets do not intersect.
Finally, $\|f - h\|_{\bigcup M,\infty}\leq \e$ because for every $x \in \bigcup M$ there exists $e \in L^1$ such that
\[ 
 x \in e , \qquad f(x) , h(x) \in \bar{B}( h(e) , \frac{\e}{3}) \quad \mbox{and} \quad \diam \bar B(h(e), \frac{\e}{3}) \leq \e .
\]

Since $h$ and $h^{-1}$ are H\"older continuous, we are able to estimate the length of each side of the $1$-dimensional complex over which $f$ is piecewise affine, in terms of a power of $\e$.

This finishes the sketch of the proof of Theorem \ref{th:skeleton}. The rest of the section is devoted to a rigorous proof of it.

Let $|\cdot|$ be the Euclidean norm in $\R^2$. All the definitions based on the norm (such as the balls and the H\"older norm) are referred to $|\cdot|$.

Let $\O \subset \R^2$ and $\a, \tilde{\a} \in (0,1]$.
Let $h \in C^{\a}( \O , \R^2)$ be a homeomorphism such that $|h|_{\a} \leq H$.
Suppose $h^{-1} \in C^{\tilde{\a}} ( h(\O) , \R^2)$ and $|h|_{\tilde{\a}} \leq \tilde{H}$.
Then the following inequalities are immediate, and will be use throughout this paper without further mention:
for all $X,Y \subset \O$,
\begin{align*}
 & \left( \frac{\dist(X,Y)}{\tilde{H}} \right)^{1/\tilde{\a}} \leq \dist (h(X) , h(Y)) \leq H \dist(X,Y) ^{\a} ,   \\
 & \left( \frac{\diam X}{\tilde{H}} \right)^{1/\tilde{\a}} \leq \diam h(X) \leq H (\diam X) ^{\a} .
\end{align*}
If $\a = \tilde{\a}=1$ and $\O$ has more than one element, then clearly $H \tilde{H} \geq 1$.

The following easy property will be useful in the proof of Theorem \ref{th:skeleton}.
\begin{lemma}\label{le:curveHolder}
Let $\O \subset \R^2$ and $\a \in (0,1]$.
Let $h \in C^{\a}( \O , \R^2)$ and $\d > 0$.
If $a, b \in \O$ satisfy
\[
 |b-a| < \left( \frac{2 \d}{|h|_{\a}} \right)^{1/\a}
\]
then $[h(a), h(b)] \subset B (h[a,b], \d)$.
\end{lemma}
\begin{proof}
If $0 \leq \l \leq 1/2$ then $|\l h(a) + (1-\l) h(b) - h (b)| < \d$,
whereas if $1/2 \leq \l \leq 1$ then $|\l h(a) + (1-\l) h(b) - h (a)| < \d$. This concludes the proof.
\end{proof}

The following lemma starts with a piecewise affine function $f$ over a segment, and constructs a piecewise affine homeomorphism that coincides with $f$ at the endpoints of the segment and whose image is contained in the image of $f$. It will be useful in the proof of Theorem \ref{th:skeleton}.
\begin{lemma}\label{le:fginjective}
Let $n \geq 1$ be a natural number, and $a_0 \neq a_n \in \R^2$.
Consider $a_1 , \ldots, a_{n-1} \in (a_0 , a_n)$ such that $a_1 < \cdots < a_{n-1}$ in the order of $(a_0 , a_n)$. Let $f : [a_0 , a_n] \to \R^2$ be a piecewise affine function over $\{ a_0, \ldots, a_n \}$. Then
\begin{enumerate}
\item\label{1} If the points $f(a_0), \ldots, f(a_n)$ are different then there exist a natural number $m \leq n$, points $b_0, \ldots, b_m \in [a_0, a_n]$ such that $a_0 = b_0 < \cdots < b_m = a_n$ and an injective function $g : [a_0 , a_n] \to \R^2$ piecewise affine over $\{ b_0, \ldots, b_m \}$ such that $g|_{\{ a_0,a_n \}} = f|_{\{ a_0,a_n \}}$ and $g [a_0 , a_n] \subset f [a_0 , a_n]$.

\item\label{2} If $f$ is injective then there exists an injective function $g : [a_0 , a_n] \to \R^2$ piecewise affine over $\{ b_0, \ldots, b_n \}$, where
\[
 b_i := a_0 + \frac{i}{n}(a_n - a_0), \qquad i= 0, \ldots, n ,
\]
such that $g|_{\{ a_0,a_n \}} = f|_{\{ a_0,a_n \}}$ and $g [a_0 , a_n] = f [a_0 , a_n]$.
\end{enumerate}
\end{lemma}
\begin{proof}
First we prove part \ref{1} by induction on $n$. When $n=1$, then $f$ is injective and we take $m=n$ and $g=f$.

Suppose $n\geq 2$. By the induction assumption applied to $f|_{[a_1,a_n]}$, there exist a natural number $m_1 \leq n-1$, points $b_0, \ldots, b_{m_1}$ with $a_1 = b_0 < \cdots < b_{m_1} = a_n$, and an injective function $g_1 : [a_1,a_n] \to \R^2$ piecewise affine over $\{ b_0, \ldots, b_{m_1} \}$ such that $g_1 |_{\{ a_1,a_n \}} = f|_{\{ a_1,a_n \}}$ and $g_1 [a_1 , a_n] \subset f [a_1 , a_n]$.
Define $g_2 : [a_0, a_n] \to \R^2$ as $g_2 |_{[a_0,a_1)} = f |_{[a_0,a_1)}$ and $g_2 |_{[a_1,a_n]} = g_1 |_{[a_1,a_n]}$.
Then $g_2$ is piecewise affine over $\{ a_0, a_1, b_1, \ldots, b_{m_1 -1},a_n \}$; moreover, $g_2 |_{\{ a_0,a_n \}} = f|_{\{ a_0,a_n \}}$ and $g_2 [a_0 , a_n] \subset f [a_0 , a_n]$. If $g_2$ is injective then we are done.

Suppose $g_2$ is not injective. Since $g_2$ is injective in $[a_0, a_1]$ and in $[a_1, a_n]$,
\[
 g_2 [a_0, a_1) \cap g_2 [a_1, a_n] \neq \varnothing .
\]
Define
\[
 c_0 := \inf \{ x \in [a_0, a_1) : g_2 (x) \in g_2 [a_1, a_n] \} .
\]
Since $g_2 [a_1, a_n]$ is closed, this $c_0 \in [a_0, a_1)$ satisfies $g_2 (c_0) \in g_2 [a_1, a_n]$.
Moreover, as $g_2$ is injective in $[a_0, a_1]$ and in $[a_1, a_n]$, there exists a unique $c_1 \in (a_1, a_n]$ such that $g_2 (c_0) = g_2 (c_1)$.
If $a_0 = c_0$ then the sets $[a_0, c_0)$ and $g_2 [a_0, c_0)$ are empty.
By definition of $c_0$,
\[
 g_2 [ a_0, c_0 ) \cap g_2 [ a_1, a_n ] = \varnothing ,
\]
and, hence, $g_2$ is injective in $[a_0, c_0) \cup [c_1, a_n]$. The function $g: [a_0 , a_n] \to \R^2$ defined by
\[
 g(x) := \left\{ \begin{array}{lcl}
  g_2(x), & \mbox{if} & x \in [a_0, c_0) \\
  g_2 ( |a_n - c_1| |a_n - c_0|^{-1} (x - c_0) + c_1 ), & \mbox{if}& x \in [c_0, a_n].
  \end{array} \right.
\]
satisfies the requirements of the statement.
This proves part \ref{1}.

Now we prove part \ref{2}. For each $k \in \{ 1, \ldots, n\}$, let $h_k : [b_{k-1}, b_k] \to [a_{k-1}, a_k]$ be the affine function such that  $h_k (b_{k-1}) = a_{k-1}$ and  $h_k (b_k) = a_k$. Define $g : [a_0 , a_n] \to \R^2$ as the only function that, for each $k \in \{ 1, \ldots, n\}$, coincides with $ f \circ h_k$ in $[b_{k-1}, b_k]$.
This $g$ satisfies the requirements of the statement, which proves part \ref{2}.
\end{proof}

The following is the main result of this section.
Given a $1$-dimensional complex $M$ and a homeomorphism $h$ defined on $\bigcup M$, we construct a piecewise affine homeomorphism over a refinement $K$ of the complex $M$ that approximates $h$ in the supremum norm.
Its proof is based on the one of Theorem 2 of Chapter 6 of Moise \cite{Moise}, but with the difference that we make an explicit construction, and estimate the lengths of the elements of $K^1$.

\begin{theorem}\label{th:skeleton}
Let $M$ be a $1$-dimensional complex in $\R^2$.
Let $\t \in (0,\pi/3]$ satisfy that $\sin \f \geq \sin \t$ for all angles $\f$ of all the triangles defined by $M$.
Let $\a , \tilde{\a} \in (0, 1]$.
Let $h \in C^{\a} (\bigcup M , \R^2)$ be a homeomorphism such
that $h^{-1} \in C^{\tilde{\a}} (h(\bigcup M) , \R^2)$.

Then there exists a constant $\e_0 >0$, depending only on
\[
 \min_{e \in M^1} |e|, \quad |h|_{\a}, \quad |h^{-1}|_{\tilde{\a}}, \quad \a, \quad \tilde{\a} ,
\]
such that for each $0 < \e \leq \e_0$
there exist a homeomorphism $f: \bigcup M \to \R^2$
and a $1$-dimensional complex $K$ such that
\[
 \bigcup K = \bigcup M , \qquad M^0 \subset K^0 , \qquad \| f-h \|_{\infty} \leq \e,
\]
$f$ is piecewise affine over $K$, coincides with $h$ in $\tilde{M}^0$,
and
\[
 B_1 \e^{b_1} < |e| < B_2 \e^{b_2} , \qquad e \in K^1 ,
\]
where
\begin{equation}\label{eq:B}
\begin{split}
 b_1 := & - \frac{2}{\a} + \frac{2}{\a^2 \tilde{\a}} + \frac{1}{\a^3 \tilde{\a}^2} , \\
 b_2 := & \frac{1}{\a} , \\
 B_1 := &
 2^{-2 - \frac{2}{\a \tilde{\a}} - \frac{1}{\a^2 \tilde{\a}} - \frac{1}{\a^2 \tilde{\a}^2}} \cdot
 3^{\frac{1}{\a} -\frac{2}{\a^2 \tilde{\a}} - \frac{1}{\a^3 \tilde{\a}^2}} \cdot
 \left[ 1 - (\sin \t)^{1/\tilde{\a}} \right]^{2/\a} \cdot
 (\sin \t)^{\frac{1}{\a \tilde{\a}} + \frac{1}{\a^2 \tilde{\a}^2} } \\
  & \cdot |h|_{\a}^{-\frac{1}{\a} - \frac{3}{\a^2 \tilde{\a}} - \frac{1}{\a^3 \tilde{\a}^2}} \cdot
  |h^{-1}|_{\tilde{\a}}^{- \frac{3}{\a \tilde{\a}} - \frac{1}{\a^2 \tilde{\a}^2} } \qquad \mbox{if } \ \a \tilde{\a} <1 , \\
 B_1:= & \frac{\sin^2 \t (1-\sin \t)}{192 |h|_{\a}^5 |h^{-1}|_{\tilde{\a}}^4} \min \left\{ 1- \sin \t , \frac{\sin^2 \t}{2 |h|_{\a} |h^{-1}|_{\tilde{\a}}} \right\}  \qquad \mbox{if } \ \a = \tilde{\a} = 1 , \\
 B_2:= & 3^{-1/\a} |h|_{\a}^{-1/\a} .
\end{split}
\end{equation}
\end{theorem}
\begin{proof}
Call $H := |h|_{\a}$ and $\tilde{H} := |h^{-1}|_{\tilde{\a}}$.
Let $0 < \e \leq \e_0$, where $\e_0$ is to be determined later.

Fix $e \in M^1$. Let $n_e$ be the only integer $x$ satisfying
\begin{equation}\label{eq:Nj}
  |e| \left( \frac{\e}{3 H} \right) ^{-1/\a} \leq x < |e| \left( \frac{\e}{3 H} \right) ^{-1/\a} +1 .
\end{equation}
Choose one of the two affine bijections from $[0,1]$ onto $e$, and equip $e$ with the total order inherited from $[0,1]$ and given by that bijection.
Now define
\[
 v_{e,i} := \inf e + \frac{i}{n_e} (\sup e - \inf e) , \qquad i= 0, \ldots , n_e .
\]
If $\e_0$ is small then
\begin{equation}\label{eq:Njgeq2}
 n_e \geq 2 .
\end{equation}

Let $L$ be the $1$-dimensional complex defined by
\[
 L:= \left\{ \{v_{e,i}\} : \ \ e \in M^1, \ \ i= 0, \ldots , n_e \right\} \cup
  \left\{ [ v_{e,i} , v_{e,i+1}] : \ \ e \in M^1, \ \ i= 0, \ldots , n_e -1 \right\} .
\]

By \eqref{eq:Nj}, \eqref{eq:Njgeq2} and the geometry of the triangles defined by $M$,
\begin{align}
 &\label{eq:|e|} \frac{1}{2} \left( \frac{\e}{3 H} \right) ^{1/\a} < |e|  \leq \left( \frac{\e}{3 H} \right) ^{1/\a} , & e \in L^1 , \\
 &\label{eq:u-v} |u-v| > \frac{\sin \t}{2} \left( \frac{\e}{3 H} \right) ^{1/\a} , & u \neq v \in \tilde{L}^0 , \\
 &\label{eq:distve} \dist (v, e) > \frac{\sin \t}{2} \left( \frac{\e}{3 H} \right) ^{1/\a} , & v \in \tilde{L}^0, \quad e\in L^1, \quad v\notin e .
\end{align}

By \eqref{eq:|e|}, for each $e \in L^1$,
\begin{equation}\label{eq:diamBhe}
 \diam h (e) \leq \frac{\e}{3} , \qquad \diam \bar{B} (h(e), \frac{\e}{3}) = \diam h(e) + \frac{2\e}{3} \leq \e .
\end{equation}

Define
\[
 \b :=
 \frac{1}{2} \left( \frac{\sin \t}{2 \tilde{H}} \right) ^{1/\tilde{\a}} \left( \frac{\e}{3H} \right) ^{\frac{1}{\a \tilde{\a}}} .
\]
If $\e_0$ is small and $\a \tilde{\a} < 1$ then
\begin{equation}\label{eq:b3leqe3}
 \b \leq \e/3 ;
\end{equation}
otherwise, if $\a = \tilde{\a} = 1$ then $H \tilde{H} \geq 1$ and equation \eqref{eq:b3leqe3} holds as well, and takes the form
\[
 \frac{\sin \t}{4 H \tilde{H}} \leq 1 .
\]

For each $v \in \tilde{L}^0$ define $N_v := \bar{B} ( h(v) , \b)$.
By \eqref{eq:u-v}, for all $u \neq v \in \tilde{L}^0$,
\[
 |h(u) - h(v)| \geq \left( \frac{|u-v|}{\tilde{H}}\right) ^{1/\tilde{\a}} > 2 \b,
\]
and, hence,
\begin{equation}\label{eq:NuNv}
 N_u \cap N_v = \varnothing .
\end{equation}
Furthermore, if $v \in \tilde{L}^0$ does not belong to $e \in L^1$, by \eqref{eq:distve} we get
\begin{equation*}
 \dist ( h(v) , h(e)) \geq \left( \frac{\dist(v,e)}{\tilde{H}} \right) ^ {1/\tilde{\a}}  > 2 \b
\end{equation*}
and, hence,
\begin{equation}\label{eq:NBempty}
 N_v \cap \bar{B} (h(e), \b) = \varnothing .
\end{equation}

Let $e \in L^1$. Equip $h(e)$ with the total order inherited from $e$ and given by the bijection $h$.
By \eqref{eq:|e|},
\begin{equation*}
 |h(\sup e) - h(\inf e)| \geq \left( \frac{|e|}{\tilde{H}} \right) ^{1/\tilde{\a}} >
 \left( \frac{1}{2 \tilde{H}} \right) ^{1/\tilde{\a}} \left( \frac{\e}{3H} \right) ^{\frac{1}{\a \tilde{\a}}} > 2 \b .
\end{equation*}
Therefore, there exist $x_e \in \p N_{\inf e} \cap h(e)$ and $y_e\in \p N_{\sup e} \cap h([h^{-1} (x_e), \sup e])$.
By \eqref{eq:NuNv}, $x_e < y_e$ in the order of $h(e)$.
In addition, by \eqref{eq:|e|},
\begin{equation}\label{eq:h-1ye-h-1xeleq}
 | h^{-1} (y_e) - h^{-1} (x_e) | < |e| \leq \left( \frac{\e}{3 H} \right) ^{1/\a}
\end{equation}
and
\begin{multline}\label{eq:h-1ye-h-1xegeq}
 | h^{-1} (y_e) - h^{-1} (x_e) | \geq \left( \frac{|y_e - x_e|}{H} \right) ^{1/\a}
 \geq  \left( \frac{| h(\sup e) - h(\inf e) | - 2 \b}{H} \right) ^{1/\a} \\
 >  \left( \frac{\left( \frac{1}{2 \tilde{H}} \right)^{1/\tilde{\a}} \left( \frac{\e}{3H} \right)^{\frac{1}{\a \tilde{\a}}} -2\b }{H} \right)^{1/\a}
  = \left( \frac{1 - (\sin \t) ^{1/\tilde{\a}}}{H} \right)^{1/\a}  \left( \frac{1}{2 \tilde{H}} \right)^{\frac{1}{\a\tilde{\a}}}  \left( \frac{\e}{3H} \right)^{\frac{1}{\a^2 \tilde{\a}}} .
\end{multline}
Define
$A_e := h ([ h^{-1} (x_e) , h^{-1} (y_e) ])$.
By the geometry of the triangles defined by $M$, for every $e \neq c \in L^1$,
\begin{equation}\label{eq:AeAc}
\begin{split}
 \dist (A_e, A_c) &\geq \left( \frac{1}{\tilde{H}} \right)^{1/\tilde{\a}} \left( \dist \left( [ h^{-1} (x_e) , h^{-1} (y_e) ] , [ h^{-1} (x_c) , h^{-1} (y_c) ] \right) \right)^{1/\tilde{\a}} \\
  &> \left( \frac{1}{\tilde{H}} \right)^{1/\tilde{\a}} \left( \sin \t \min \left\{ \displaystyle \min_{d \in L^1} | \sup d - h^{-1} (y_d)| , \min_{d \in L^1} | h^{-1} (x_d) - \inf d| \right\} \right)^{1/\tilde{\a}} \\
  &\geq \left( \frac{1}{\tilde{H}} \right)^{1/\tilde{\a}} \left( \sin \t \left[ \displaystyle \frac{1}{H} \min \left\{\min_{d \in L^1} | h(\sup d) - y_d| , \min_{d \in L^1} | x_d - h(\inf d)| \right\} \right]^{1/\a} \right)^{1/\tilde{\a}} \\
  &= \left( \frac{\sin \t}{\tilde{H}} \right) ^{1/\tilde{\a}} \left( \frac{\b}{H} \right) ^{\frac{1}{\a \tilde{\a}}} .
\end{split}
\end{equation}

Define
\begin{equation*}
  \d := \frac{1}{2} \left( \frac{\sin \t}{\tilde{H}} \right) ^{1/\tilde{\a}} \left( \frac{\b}{H} \right) ^{\frac{1}{\a \tilde{\a}}}
  = \frac{1}{2} \left( \frac{\sin \t}{\tilde{H}} \right) ^{1/\tilde{\a}} \left( \frac{1}{2H} \right) ^{\frac{1}{\a \tilde{\a}}}
  \left( \frac{\sin \t}{2 \tilde{H}} \right) ^{\frac{1}{\a\tilde{\a}^2}} \left( \frac{\e}{3H} \right) ^{\frac{1}{\a^2 \tilde{\a}^2}} .
\end{equation*}
By \eqref{eq:AeAc},
\begin{equation}\label{eq:BAeBAc}
\bar{B} (A_e, \d) \cap \bar{B} (A_c, \d) = \varnothing, \qquad e \neq c \in L^1 .
\end{equation}
For the reader's convenience, from now on in this proof we will write how the estimates behave with $\e > 0$. Precisely, if a quantity $\varrho$ has been defined, and $a \in \R$, when we write $\varrho \sim \e^a$ we mean that there exists a constant $c >0$ depending only on $\t, H, \tilde{H}, \a, \tilde{\a}$ such that $\varrho = c \e^a$. The constant $c$ will have been calculated explicitly, so that this notation is only a reminder.
For example, we have already showed that 
\[
 \b \sim \e^{\frac{1}{\a \tilde{\a}}}, \qquad \d \sim \e^{\frac{1}{\a^2 \tilde{\a}^2}} .
\]
Thus, if $\e_0$ is small and $\a \tilde{\a} < 1$ then
\begin{equation}\label{eq:dleqb3}
 \d \leq \b  ;
\end{equation}
otherwise, if $\a = \tilde{\a} =1$ then $H \tilde{H} \geq 1$ 
and equation \eqref{eq:dleqb3} holds as well, and takes the form
\[
 \frac{\sin \t}{2 H \tilde{H}} \leq 1 .
\]

Let $N$ be the only integer $x$ satisfying
\begin{equation}\label{eq:Ne}
 \left( \frac{\e}{3 H} \right) ^{1/\a} \left( \frac{2 \d}{H} \right)^{-1/\a} \leq x < \left( \frac{\e}{3 H} \right) ^{1/\a} \left( \frac{2 \d}{H} \right)^{-1/\a} +1 .
\end{equation}
Note that 
\[
 \left( \frac{\e}{3 H} \right) ^{1/\a} \left( \frac{2 \d}{H} \right)^{-1/\a} \sim \e^{\frac{1}{\a} - \frac{1}{\a^3 \tilde{\a}^2}} .
\]
Thus, if $\e_0$ is small and $\a \tilde{\a} < 1$ then
\begin{equation}\label{eq:g>1}
 \left( \frac{\e}{3 H} \right) ^{1/\a} \left( \frac{2 \d}{H} \right) ^{-1/\a} > 1 ;
\end{equation}
otherwise, if $\a = \tilde{\a} = 1$ then $H \tilde{H} \geq 1$ and equation
\eqref{eq:g>1} holds as well, and takes the form
\[
  \frac{2 H^2 \tilde{H}^2}{\sin^2 \t} > 1 .
\]
Inequalities \eqref{eq:Ne} and \eqref{eq:g>1} show that
\begin{equation}\label{eq:Negeq2}
 N \geq 2 \quad \mbox{and} \quad \frac{1}{N} > \frac{1}{2} \left( \frac{\e}{3 H} \right) ^{-1/\a} \left( \frac{2 \d}{H} \right) ^{1/\a} .
\end{equation}
Define
\begin{equation}\label{eq:wek}
 w_{e,k} := h^{-1} (x_e) + \frac{k}{N} [ h^{-1} (y_e) - h^{-1} (x_e) ] , \qquad e \in L^1, \quad k= 0, \ldots, N .
\end{equation}

Fix $e \in L^1$ and $k \in \{ 0, \ldots , N - 1 \}$.
By \eqref{eq:h-1ye-h-1xeleq} and \eqref{eq:Ne},
\begin{equation}\label{eq:lengthw}
 | w_{e,k+1} - w_{e,k}| < \left( \frac{2 \d}{H} \right) ^{1/\a} ,
\end{equation}
By \eqref{eq:Negeq2} and \eqref{eq:h-1ye-h-1xegeq},
\begin{equation}\label{eq:wleq}
\begin{split}
 |w_{e,k+1} - w_{e,k}| & = \frac{1}{N} |h^{-1}(y_e) - h^{-1}(x_e)| \\
 & > \frac{1}{2} \left( \frac{\e}{3H} \right)^{-\frac{1}{\a} + \frac{1}{\a^2 \tilde{\a}}} \left( \frac{2 \d}{H} \right)^{1/\a}
 \left( \frac{1 - (\sin \t)^{1/\tilde{\a}}}{H} \right)^{1/\a} \left( \frac{1}{2 \tilde{H}} \right)^{\frac{1}{\a \tilde{\a}}}.
\end{split}
\end{equation}
By Lemma \ref{le:curveHolder} and \eqref{eq:lengthw},
\begin{equation}\label{eq:hsubsetB}
 [h(w_{e,k}), h(w_{e,k+1})] \subset B (h[w_{e,k},w_{e,k+1}], \d) .
\end{equation}

Fix $e \in L^1$. Define $g_e : [h^{-1}(x_e), h^{-1}(y_e)] \to \R^2$ as the piecewise affine function over $\{ w_{e,0}, \ldots, w_{e,N} \}$ that coincides with $h$ in $\{ w_{e,0}, \ldots, w_{e,N} \}$.
Define
\begin{align*}
 & p_e := \sup \left\{ x \in [h^{-1}(x_e), h^{-1}(y_e)] : g_e (x) \in \p N_{\inf e} \right\} , \\
 & q_e := \inf \left\{ x \in [p_e, h^{-1}(y_e)] : g_e (x) \in \p N_{\sup e} \right\} .
\end{align*}
The point $p_e$ is well-defined since $x_e \in \p N_{\inf e}$,
and the point $q_e$ is also well-defined since $y_e \in \p N_{\sup e}$.
Moreover, $g_e (p_e) \in \p N_{\inf e}$ since $\p N_{\inf e}$ is closed, and
$g_e (q_e) \in \p N_{\sup e}$ since $\p N_{\sup e}$ is closed.
By \eqref{eq:NuNv}, $p_e < q_e$ in the order of $e$.
Furthermore, by definition of $p_e$, $q_e$ and the continuity of $g_e$,
\begin{equation}\label{eq:gexnotin}
 g_e (x) \notin N_{\inf e} \cup N_{\sup e}, \qquad x \in (p_e,q_e) .
\end{equation}
In addition, by \eqref{eq:hsubsetB},
\begin{equation}\label{eq:fsubsetB}
\begin{split}
 g_e [ p_e , q_e ] & \subset g_e [h^{-1}(x_e), h^{-1}(y_e)]
 = g_e \left( \bigcup_{k=0}^{N -1} [ w_{e,k} , w_{e,k+1} ] \right) \\
 & = \bigcup_{k=0}^{N -1} [ h(w_{e,k}) , h(w_{e,k+1}) ]
 \subset \bigcup_{k=0}^{N -1} B ( h [ w_{e,k} , w_{e,k+1} ] , \d ) = B ( A_e , \d ) .
\end{split}
\end{equation}
Let $m_e \in \N$ and $u_{e,0}, \ldots, u_{e,m_e}$ be such that
\[
 [p_e, q_e] \cap \{ w_{e,0}, \ldots, w_{e,N} \} =  \{ u_{e,0}, \ldots, u_{e,m_e} \}
\]
and $u_{e,0} < \cdots < u_{e,m_e}$. Of course,
\[
 1 \leq m_e \leq N , \quad u_{e,0} = p_e, \quad u_{e,m_e} = q_e .
\]
By \eqref{eq:gexnotin},
\begin{align*}
 & g_e (u_{e,0}) \in \p N_{\inf e}, \qquad g_e (u_{e,m_e}) \in \p N_{\sup e} , \\
 & g_e (u_{e,k}) = h (u_{e,k}) \notin \p N_{\inf e} \cup \p N_{\sup e} , \qquad k = 1, \ldots, m_e -1.
\end{align*}
This, \eqref{eq:NuNv} and the injectivity of $h$ imply that the points $g_e (u_{e,0}) , \ldots, g_e (u_{e,m_e})$
are different.
Therefore, by Lemma \ref{le:fginjective} applied to $g_e |_{[p_e,q_e]}$, there exist an injective function $f_e : [p_e, q_e] \to \R^2$ and a natural number $M_e$ with
\begin{equation}\label{eq:Me}
 1 \leq M_e \leq m_e \leq N
\end{equation}
such that
\begin{equation}\label{eq:fe}
 f_e (p_e) = g_e (p_e), \quad  f_e (q_e) = g_e(q_e), \quad f_e [p_e , q_e] \subset g_e [p_e , q_e] ,
\end{equation}
and $f_e$ is piecewise affine over $\{ b_{e,0}, \ldots, b_{e,M_e} \}$, where
\[
  b_{e,k} := p_e + \frac{k}{M_e} ( q_e - p_e ) , \qquad k= 0, \ldots, M_e .
\]
This construction is represented in Figure \ref{fig:loop}, in the particular case where $m_e = 4$ and $M_e = 3$.
\begin{figure}
\begin{center}
\input{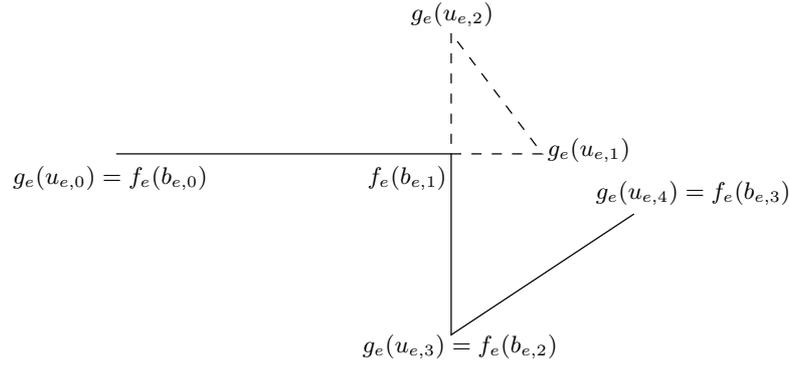}
\caption{\label{fig:loop}$g_e[p_e,q_e]$ is represented by dashed lines; $f_e[p_e,q_e]$ is represented by solid lines.}
\end{center}
\end{figure}
Note also that by \eqref{eq:|e|},
\begin{equation}\label{eq:geqe-gepeleq}
 |g_e (q_e) - g_e (p_e)| \geq |h(\sup e) - h(\inf e)| - 2 \b >
  \left( \frac{1}{2 \tilde{H}} \right)^{1/\tilde{\a}} \left( {\frac{\e}{3H}} \right)^{\frac{1}{\a \tilde{\a}}} \left[ 1- (\sin \t)^{1/\tilde{\a}} \right] .
\end{equation}
In addition, by \eqref{eq:h-1ye-h-1xeleq},
\begin{equation}\label{eq:qe-peleq}
 |q_e - p_e| \leq |h^{-1}(y_e) - h^{-1}(x_e)| < \left( \frac{\e}{3H} \right)^{1/\a} .
\end{equation}
Now we estimate $|q_e - p_e|$ from below. The analysis of this estimate varies according to whether $\a = \tilde{\a}=1$ or $\a \tilde{\a}<1$.
Suppose, first, that $\a = \tilde{\a}=1$.
We distinguish three cases according to the relative position of the points $p_e$ and $q_e$.

\begin{enumerate}
\item\label{case1pq} $p_e, q_e \in [w_{e,j}, w_{e,j+1}]$ for some $0\leq j\leq N -1$.

\item\label{case2pq} $p_e \in [w_{e,j}, w_{e,j+1}]$ and $q_e \in [w_{e,k}, w_{e,k+1}]$ for some $0\leq j,k\leq N -1$ with $j+1 < k$.

\item\label{case3pq} $p_e \in [w_{e,k-1}, w_{e,k}]$ and $q_e \in [w_{e,k}, w_{e,k+1}]$ for some $0\leq k\leq N -1$.
\end{enumerate}

Suppose Case \ref{case1pq}. Then there exist numbers $\l,\mu$ such that $0\leq \l < \mu \leq 1$ and
\begin{equation}\label{eq:case1peqe}
\begin{split}
  p_e = (1-\l) w_{e,j} + \l w_{e,j+1}, & \qquad g_e (p_e) = (1-\l) h(w_{e,j}) + \l h(w_{e,j+1}) , \\
  q_e = (1-\mu) w_{e,j} + \mu w_{e,j+1}, & \qquad g_e (q_e) = (1-\mu) h(w_{e,j}) + \mu h(w_{e,j+1}) .
 \end{split}
\end{equation}
Then
\begin{equation}\label{eq:case1pe-qe}
 |q_e -p_e| = (\mu - \l) |w_{e,j+1} - w_{e,j}| , \quad
 |g_e(q_e) - g_e(p_e)| = (\mu - \l) |h(w_{e,j+1}) - h(w_{e,j})| .
\end{equation}
By \eqref{eq:lengthw},
\begin{equation}\label{eq:hwleqLip}
 |h(w_{e,j+1}) - h(w_{e,j})| \leq H |w_{e,j+1} - w_{e,j}| < 2 \d .
\end{equation}
Therefore, by \eqref{eq:geqe-gepeleq} and \eqref{eq:hwleqLip},
\begin{equation}\label{eq:l-mLip}
 \mu - \l = \frac{|g_e(q_e) - g_e(p_e)|}{|h(w_{e,j+1}) - h(w_{e,j})|}
 > \left( \frac{1}{2 \tilde{H}} \right) \left( {\frac{\e}{3H}} \right) \left[ 1- \sin \t \right] \frac{1}{2 \d} = 2 H \tilde{H} \frac{1 - \sin \t}{\sin^2 \t}.
\end{equation}
So inequalities \eqref{eq:wleq} and \eqref{eq:l-mLip} imply
\begin{equation}\label{eq:qe-pe1Lip}
 |q_e - p_e| > \frac{(1 - \sin \t)^2}{24 H^3 \tilde{H}^2} \e .
\end{equation}

Suppose Case \ref{case2pq}.
Then, by \eqref{eq:wleq},
\begin{equation*}
  | q_e - p_e | \geq  | w_{e,k} - w_{e,j+1} | \geq  | w_{e,k} - w_{e,k-1} | >
 \frac{\d (1 - \sin \t)}{2 H^2 \tilde{H}} .
\end{equation*}

Suppose Case \ref{case3pq}. Then there exist $\l,\mu \in [0,1]$ such that
\begin{equation}\label{eq:case3peqe}
\begin{split}
  p_e = (1-\l) w_{e,k-1} + \l w_{e,k}, & \qquad  g_e (p_e) = (1-\l) h(w_{e,k-1}) + \l h(w_{e,k}) , \\
  q_e = (1-\mu) w_{e,k} + \mu w_{e,k+1}, & \qquad  g_e (q_e) = (1-\mu) h(w_{e,k}) + \mu h(w_{e,k+1}) .
\end{split}
\end{equation}
Then
\begin{equation}\label{eq:case3pe-qe}
 |q_e -p_e| = \left|(1-\mu - \l) w_{e,k} + \mu w_{e,k+1} + (\l-1) w_{e,k-1} \right| =
 (1 + \mu - \l) \frac{| h^{-1}(y_e) - h^{-1}(x_e) |}{N} .
\end{equation}
On the other hand, by \eqref{eq:lengthw},
\begin{equation}\label{eq:geq_e-g_ep_e3Lip}
\begin{split}
 |g_e(q_e) - g_e(p_e)| & =  \left| (1-\l) [h(w_{e,k})-h(w_{e,k-1})] + \mu [h(w_{e,k+1})-h(w_{e,k})] \right| \\
 & \leq (1-\l) H |w_{e,k}-w_{e,k-1}| + \mu H |w_{e,k+1}-w_{e,k}| < 2 (1-\l+\mu) \d.
\end{split}
\end{equation}
Therefore, inequalities \eqref{eq:geq_e-g_ep_e3Lip} and \eqref{eq:geqe-gepeleq} show that
\begin{equation}\label{eq:1-l+mLip}
 1 - \l + \mu > \frac{|g_e(q_e)-g_e(p_e)|}{2 \d} >
 \frac{1}{2 \tilde{H}} {\frac{\e}{3H}} \left( 1- \sin \t \right) \frac{1}{2 \d} =  2 H \tilde{H} \frac{1 - \sin \t}{\sin^2 \t}.
\end{equation}
So inequalities \eqref{eq:wleq} and \eqref{eq:1-l+mLip} imply \eqref{eq:qe-pe1Lip}.

The analysis of the three cases above shows that when $\a = \tilde{\a} = 1$ then
\begin{equation}\label{eq:caselip}
 |q_e -p_e| \geq \frac{1 - \sin \t}{24 H^3 \tilde{H}^2} \e \min \{ 1 - \sin \t , \frac{\sin^2 \t}{2 H \tilde{H}} \} .
\end{equation}

Now we estimate $|q_e - p_e|$ from below in the case $\a \tilde{\a} < 1$. We distinguish the same three cases as above.

Suppose Case \ref{case1pq}. Then there exist numbers $\l,\mu$ such that $0\leq \l < \mu \leq 1$ and identities \eqref{eq:case1peqe} and \eqref{eq:case1pe-qe} hold.
By \eqref{eq:lengthw},
\begin{equation}\label{eq:hwleq}
 |h(w_{e,j+1}) - h(w_{e,j})| \leq H |w_{e,j+1} - w_{e,j}|^{\a} < 2 \d .
\end{equation}
Therefore, by \eqref{eq:geqe-gepeleq} and \eqref{eq:hwleq},
\begin{equation*}
 \mu - \l = \frac{|g_e(q_e) - g_e(p_e)|}{|h(w_{e,j+1}) - h(w_{e,j})|}
 > \left( \frac{1}{2 \tilde{H}} \right)^{1/\tilde{\a}} \left( {\frac{\e}{3H}} \right)^{\frac{1}{\a \tilde{\a}}} \left[ 1- (\sin \t)^{1/\tilde{\a}} \right] \frac{1}{2 \d} \sim \e^{\frac{1}{\a \tilde{\a}} - \frac{1}{\a^2 \tilde{\a}^2}} ,
\end{equation*}
which is a contradiction if $\e_0$ is small, since this would imply $\mu - \l > 1$.

Suppose Case \ref{case2pq}. Then there exist $\l,\mu \in [0,1]$ such that
\begin{equation*}
 g_e (p_e) = (1-\l) h(w_{e,j}) + \l h(w_{e,j+1}) , \qquad
 g_e (q_e) = (1-\mu) h(w_{e,k}) + \mu h(w_{e,k+1}) .
\end{equation*}
By \eqref{eq:lengthw},
\begin{equation}\label{eq:case2gepegeqe}
\begin{split}
 |g_e (q_e) - g_e (p_e)| & \leq H \left( |w_{e,k} - w_{e,j}|^{\a} + \mu |w_{e,k+1} - w_{e,k}|^{\a} + \l |w_{e,j+1} - w_{e,j}|^{\a} \right) \\
 & \leq 2 \d \left[ (k-j)^{\a} + \mu + \l \right]  \leq 2 \d \left[ (k-j)^{\a} + 2 \right].
\end{split}
\end{equation}
Now define
\begin{align*}
 & \g_1 := \left( \frac{1}{2 \tilde{H}} \right)^{1/\tilde{\a}} \left( {\frac{\e}{3H}} \right)^{\frac{1}{\a \tilde{\a}}} \left[ 1- (\sin \t)^{1/\tilde{\a}} \right] \frac{1}{2\d} \sim \e^{\frac{1}{\a \tilde{\a}} - \frac{1}{\a^2 \tilde{\a}^2}}, \\
 & \g_2 :=  \frac{1}{2} \left( \frac{\e}{3H} \right)^{-\frac{1}{\a} + \frac{1}{\a^2 \tilde{\a}}} \left( \frac{2 \d}{H} \right)^{1/\a}
 \left( \frac{1 - (\sin \t)^{1/\tilde{\a}}}{H} \right)^{1/\a} \left( \frac{1}{2 \tilde{H}} \right)^{\frac{1}{\a \tilde{\a}}} \sim \e^{-\frac{1}{\a} + \frac{1}{\a^2 \tilde{\a}} + \frac{1}{\a^3 \tilde{\a}^2}}.
\end{align*}
Inequalities \eqref{eq:case2gepegeqe} and \eqref{eq:geqe-gepeleq} show that
\begin{equation}\label{eq:k-j}
  (k-j)^{\a} + 2 \geq \frac{|g_e (q_e) - g_e (p_e)|}{2 \d} \geq \g_1 .
\end{equation}
Now, by \eqref{eq:wleq}
\begin{equation}\label{eq:qeperefine}
 |q_e - p_e| \geq |w_{e,k} - w_{e,j+1}| \geq (k-j-1) \g_2 .
\end{equation}
As $\e_0$ is small, then
$\g_1 \geq \max \{ 4, \left( 2^{-1/\a} + 3^{-1/\a} \right)^{-\a} \}$
and, hence,
\begin{equation}\label{eq:gamma12}
 \g_1 - 2 \geq \frac{\g_1}{2} \qquad \mbox{and} \qquad \g_2 \left( \frac{\g_1}{2} \right)^{1/\a} - \g_2 \geq \g_2 \left( \frac{\g_1}{3} \right)^{1/\a} .
\end{equation}
In total, equations \eqref{eq:k-j}, \eqref{eq:qeperefine} and \eqref{eq:gamma12} imply
\begin{equation}\label{eq:qe-pe2}
 |q_e - p_e| \geq \g_2 \left( \frac{\g_1}{3} \right)^{1/\a} .
\end{equation}

Suppose Case \ref{case3pq}. Then there exist $\l,\mu \in [0,1]$ such
that identities \eqref{eq:case3peqe} and \eqref{eq:case3pe-qe} hold.
On the other hand, by \eqref{eq:lengthw},
\begin{equation}\label{eq:geq_e-g_ep_e3}
\begin{split}
 |g_e(q_e) - g_e(p_e)| &=  \left| (1-\l) [h(w_{e,k})-h(w_{e,k-1})] + \mu [h(w_{e,k+1})-h(w_{e,k})] \right| \\
 &\leq (1-\l) H |w_{e,k}-w_{e,k-1}|^{\a} + \mu H |w_{e,k+1}-w_{e,k}|^{\a} < 2 (1-\l+\mu) \d.
\end{split}
\end{equation}
Therefore, inequalities \eqref{eq:geq_e-g_ep_e3} and \eqref{eq:geqe-gepeleq} show that
\begin{equation*}
 1 - \l + \mu > \frac{|g_e(q_e)-g_e(p_e)|}{2 \d} >
\left( \frac{1}{2 \tilde{H}} \right)^{1/\tilde{\a}} \left( {\frac{\e}{3H}} \right)^{\frac{1}{\a \tilde{\a}}} \left[ 1- (\sin \t)^{1/\tilde{\a}} \right] \frac{1}{2 \d} \sim \e^{\frac{1}{\a \tilde{\a}} - \frac{1}{\a^2 \tilde{\a}^2}},
\end{equation*}
which is a contradiction if $\e_0$ is small, since this would imply $1 - \l + \mu > 2$.

Define
\[
 \tau := \left\{
     \begin{array}{lll}
 \frac{1}{2} \left( \frac{[1-(\sin \t)^{1/\tilde{\a}}]^2}{3H^2} \right) ^{1/\a}
 \left( \frac{1}{2\tilde{H}} \right) ^{\frac{2}{\a \tilde{\a}}}
 \left( \frac{\e}{3H} \right) ^{-\frac{1}{\a}+\frac{2}{\a^2 \tilde{\a}}} & \mbox{if} & \a \tilde{\a} < 1,\\
 \min \left\{ \frac{(1-\sin\t)^2}{24 H^3 \tilde{H}^2} , \frac{\sin^2 \t (1-\sin\t)}{48 H^4 \tilde{H}^3} \right\} \e &  \mbox{if} & \a = \tilde{\a} = 1.
     \end{array} \right.
\]
From the analysis of the three cases above, in particular from inequalities \eqref{eq:caselip} and \eqref{eq:qe-pe2}, we conclude that, regardless of the value of $\a$ and $\tilde{\a}$,
\begin{equation}\label{eq:>tau}
 |q_e - p_e| > \tau , \qquad e \in L^1.
\end{equation}

Let $K$ be the $1$-dimensional complex whose set of vertices is
\[
 \tilde{K}^0 = \left\{ b_{e,k} : \quad e \in L^1, \quad k= 0, \ldots , M_e \right\} \cup \tilde{L}^0
\]
and whose set of edges is
\begin{align*}
  K^1 = & \left\{ [ b_{e,k} , b_{e,k+1} ] : \quad e \in L^1, \quad k= 0, \ldots , M_e -1 \right\} \\
 & \cup \left\{ [ \inf e , b_{e,0} ] : e \in L^1 \right\} \cup \left\{ [ b_{e,M_e} , \sup e ] : e \in L^1 \right\}.
\end{align*}
Let $e \in M^1$. Then by \eqref{eq:qe-peleq}, \eqref{eq:Me}, \eqref{eq:>tau} and \eqref{eq:Negeq2},
\begin{equation}\label{eq:lengthb}
 \frac{\tau}{2} \left( \frac{\e}{3 H} \right) ^{-1/\a} \left( \frac{2 \d}{H} \right) ^{1/\a} < | b_{e,k+1} - b_{e,k}| < \left( \frac{\e}{3H} \right) ^{1/\a} , \qquad e \in L^1, \quad k= 0, \ldots , M_e - 1.
\end{equation}
Fix $e \in L^1$. By \eqref{eq:|e|},
\begin{equation}\label{eq:b-inf}
 \max\{ | b_{e,0} - \inf e | , | \sup e - b_{e,M_e} | \} < |e| \leq \left( \frac{\e}{3H} \right) ^{1/\a} .
\end{equation}
Moreover,
\begin{equation}\label{eq:b-inf>}
\begin{split}
 &  | b_{e,0} - \inf e | \geq |h^{-1}(x_e) - \inf e| \geq \left( \frac{x_e - h(\inf e)}{H} \right) ^{1/\a} =  \left( \frac{\b}{H} \right) ^{1/\a}, \\
 &  | \sup e - b_{e,M_e} | \geq |\sup e - h^{-1}(y_e)| \geq \left( \frac{h(\sup e) - y_e}{H} \right) ^{1/\a} =  \left( \frac{\b}{H} \right) ^{1/\a} .
\end{split}
\end{equation}

If $\e_0$ is small and $\a \tilde{\a} < 1$ then
\begin{equation}\label{eq:b42}
 \e^{- \frac{2}{\a} + \frac{2}{\a^2 \tilde{\a}} + \frac{1}{\a^3 \tilde{\a}^2}} \sim \frac{\tau}{2} \left( \frac{\e}{3 H} \right) ^{-1/\a} \left( \frac{2 \d}{H} \right) ^{1/\a} < \left( \frac{\b}{H} \right)^{1/\a} \sim \e^{\frac{1}{\a^2 \tilde{\a}}};
\end{equation}
otherwise, if $\a = \tilde{\a} =1$ then $H \tilde{H} \geq 1$ and equation \eqref{eq:b42} holds as well, and takes the form
\[
 \frac{(1-\sin \t)\sin \t}{16 H^3 \tilde{H}^3} \min \left\{ 1 - \sin \t , \frac{\sin^2 \t}{2 H \tilde{H}} \right\} < 1 .
\]
In total, by \eqref{eq:lengthb}, \eqref{eq:b-inf}, \eqref{eq:b-inf>} and \eqref{eq:b42},
\[
 \frac{\tau}{2} \left( \frac{\e}{3 H} \right) ^{-1/\a} \left( \frac{2 \d}{H} \right) ^{1/\a} < |e| < \left( \frac{\e}{3H} \right) ^{1/\a} , \qquad e \in K^1 .
\]
Note that
\[
 \frac{\tau}{2} \left( \frac{\e}{3 H} \right) ^{-1/\a} \left( \frac{2 \d}{H} \right) ^{1/\a}
 = B_1 \e^{b_1} ,
\]
where $B_1$ and $b_1$ are defined in \eqref{eq:B}.

Define $f : \bigcup M \to \R^2$ as the only function that, for each $e \in L^1$,
\begin{align*}
 & f|_{[p_e,q_e]} = f_e , \\
 & f|_{[\inf e,p_e]} \quad \mbox{is affine with} \quad f(\inf e) = h(\inf e), \quad f(p_e) = g_e(p_e) , \\
 & f|_{[q_e,\sup e]} \quad \mbox{is affine with} \quad f(q_e) = g_e(q_e), \quad f(\sup e) = h(\sup e) .
\end{align*}
Then $f$ is piecewise affine over $K$.
We now prove that $f$ is a homeomorphism. Since $f$ is continuous and $\bigcup M$ is compact, it suffices to show that $f$ is injective.
The proof of the injectivity of $f$ involves the consideration of several cases.

We have already seen that $f$ is injective in $[ p_e , q_e ]$
for each $e \in L^1$, since so is $f_e$.

Now consider $e, c \in L^1$ such that $\sup e = \inf c$. Then $f$ is piecewise affine over $\{ q_e , \sup e , p_c \}$,
\begin{equation*}
 f(q_e) = y_e , \quad f(\sup e) = h(\sup e) , \quad f(p_c) = x_c , \quad y_e, x_c \in \p N_{\sup e}, \quad y_e \neq x_c ;
\end{equation*}
the last inequality is due to \eqref{eq:fsubsetB} and \eqref{eq:BAeBAc}.
Thanks to the geometry of the Euclidean ball, this implies that $f$ is injective in $[q_e , \sup e] \cup [\inf c, p_c]$.

By \eqref{eq:fe} and \eqref{eq:fsubsetB}, for every $e \in L^1$,
\begin{equation}\label{eq:fsubsetB2}
 f [ p_e , q_e ] =  f_e [ p_e , q_e ] \subset  g_e [ p_e , q_e ] \subset B ( A_e , \d ) .
\end{equation}
Therefore, \eqref{eq:fsubsetB2} and \eqref{eq:BAeBAc} demonstrate that
\[
  f [ p_e , q_e ] \cap  f [ p_c , q_c ] = \varnothing , \qquad e \neq c \in L^1 .
\]

If $e,c \in L^1$ satisfy $\sup e = \inf c$ then
\begin{equation}\label{eq:fsubsetN}
 f ( [ q_e, \sup e ] \cup [ \inf c , p_c ]  )
 = [ g_e(q_e), h(\sup e) ] \cup [ h(\inf c) , g_c(p_c) ]
 \subset N_{\sup e} .
\end{equation}
Therefore, \eqref{eq:NuNv} demonstrates that
\[
 f ( [q_e , \sup e] \cup [\inf c, p_c] ) \cap f ( [q_{e'} , \sup e'] \cup [\inf c', g_{c'}(p_{c'})] ) =\varnothing.
\]
whenever $e,c , e',c' \in L^1$ satisfy
\[
 \sup e = \inf c , \quad \sup e' = \inf c', \quad (e,c) \neq (e',c') .
\]

Suppose finally that there exist
\[
 e, c, d \in L^1, \quad x \in [q_e , \sup e] \cup [\inf c, p_c], \quad y \in [ p_d , q_d ]
\]
such that $f(x) = f(y)$ and $\sup e = \inf c$. By \eqref{eq:fsubsetN},
\begin{equation}\label{eq:fxin}
f(x) \in N_{\sup e} .
\end{equation}
By \eqref{eq:fsubsetB2} and \eqref{eq:dleqb3},
\begin{equation}\label{eq:fyinBB}
 f(y) \in B (A_d, \d) \subset \bar{B} (h(d), \d) \subset \bar{B} (h(d), \b)  .
\end{equation}
Equations \eqref{eq:NBempty}, \eqref{eq:fxin} and \eqref{eq:fyinBB} imply $d \in \{c,e\}$.
If $d=e$ then equations \eqref{eq:gexnotin} and \eqref{eq:NuNv} show that $y = q_e$.
If $d=c$ then equations \eqref{eq:gexnotin} and \eqref{eq:NuNv} show that $y = p_c$.
In either case, $d \in \{c,e\}$ implies $y \in \{ q_e , p_c \}$.
Since $f$ is injective in $[q_e , \sup e] \cup [\inf c, p_c]$ then $x=y$.
This proves that $f$ is a homeomorphism.

Finally, we prove that $\| f-h \|_{\infty} \leq \e$.
Take $e \in L^1$ and suppose $x \in [ \inf e , p_e] \cup [ q_e , \sup e ]$. Then
\[
 f(x) \in [ h(\inf e) , g_e(p_e) ] \cup [ g_e(q_e) , h(\sup e) ] ,
\]
and, hence, by \eqref{eq:b3leqe3}, $\dist ( f(x) , h(e) ) \leq \b \leq \e/3$.
Now suppose $x \in [ p_e , q_e ]$. By \eqref{eq:fsubsetB2},
\[
 f(x) \in B(A_e, \d) \subset B(h(e),\d),
\]
and, hence, by \eqref{eq:dleqb3} and \eqref{eq:b3leqe3},
$\dist ( f(x) , h(e) ) < \d \leq \e/3$.
Therefore, for each $x \in \bigcup M$ there exists $e \in L^1$ such that $x \in e$ and, hence,
\[
f(x) , h(x) \in \bar{B}( h(e) , \e/3) ;
\]
consequently, by \eqref{eq:diamBhe}, $| f(x) - h(x) | \leq \e$. This concludes the proof.
\end{proof}

\section{Extension from the boundary to the whole triangle}\label{se:extension}

In this section we construct a piecewise affine homeomorphism on a triangle that extends a given
piecewise affine homeomorphism defined on the boundary of the triangle. The existence of such an extension is a consequence of the so-called \emph{piecewise linear} (or \emph{piecewise affine)} \emph{Schoenflies theorem} (e.g., Theorem III.1.C of Bing \cite{Bing}), but in this section we construct a triangulation explicitly and measure the angles and the side lengths of that triangulation.
The construction that we present in this section is based on an idea by Gupta and Wenger \cite{GW} that we explain in the following paragraphs.

Let $\D$ be a triangle.
Let $h: \p\D \to \R^2$ be a piecewise affine homeomorphism over a $1$-dimensional complex $M$ such that $\bigcup M = \p \D$.
Let $a_0, \ldots, a_{w-1} \in \R^2$ satisfy $\tilde{M}^0 = \{ a_0, \ldots, a_{w-1} \}$ and $a_0 \cdots a_{w-1}$ is a well-defined polygon, which in fact coincides with $\D$.
Since $h$ is a piecewise affine homeomorphism, then $Q:= h(a_0) \cdots h(a_{w-1})$ is a well-defined polygon whose boundary coincides with $h(\p\D)$.
Note that $Q$ need not be a well-defined $w$-gon; indeed, this happens precisely when $h$ is piecewise affine over another $1$-dimensional complex $N$ such that $\Card N^0 < \Card M^0$.

It may well happen that for some $i,j,k \in \{0,\ldots,w-1\}$, the triangle $h(a_i) h(a_j) h(a_k)$ is a well-defined triangle contained in $Q$, but $a_i a_j a_k$ is not a well-defined triangle.
This occurs precisely when the points $a_i , a_j , a_k$ are aligned.
But a small perturbation of the vertices of $M$ solves this problem.
Indeed, we can define points $b_0, \ldots, b_{w-1} \in \R^2$ such that $b_i$ is close to $a_i$ for all $i \in \{0, \ldots , w-1\}$, and $\D' := b_0 \cdots b_{w-1}$ is a well-defined $w$-gon that, in addition, is convex and contained in the interior of $\D$.
Accordingly, we perturb the polygon $Q$, that is, we choose points $x_0, \ldots, x_{w-1} \in \R^2$ such that $x_i$ is close to $h(a_i)$ for all $i \in \{0, \ldots , w-1\}$, and $Q' := x_0 \cdots x_{w-1}$ is a well-defined $w$-gon.
Now, if for some $i,j,k \in \{0,\ldots,w-1\}$, the triangle $x_i x_j x_k$ is a well-defined triangle contained in $Q'$ then $b_i b_j b_k$ is a well-defined triangle triangle contained in $\D'$; this is because $\D'$ is a well-defined convex $w$-gon.
This property can be rephrased as: every triangulation in $Q'$ induces a triangulation in $\D'$ (these triangulations are sometimes called \emph{isomorphic}).
Next, we triangulate $Q'$ without adding any extra vertex, that is to say, there exists a triangulation of $Q'$ whose set of vertices is $\{b_0, \ldots, b_{w-1}\}$; the existence of such a triangulation is a well-known result.
As explained before, this triangulation of $Q'$ induces a triangulation of $\D'$.

Now we have to triangulate $Q \setminus \inte{Q'}$ in such a way that this triangulation induces a triangulation in $\D \setminus \inte{\D'}$.
But this is immediate if the perturbation of $\D$ and $Q$ has been done in such a way that $a_k a_{k+1} b_{k+1} b_k$ and $h(a_k) h(a_{k+1}) x_{k+1} x_k$ are well-defined convex quadrilaterals, for all $k \in \{0, \ldots, w-1\}$ (here we have defined $a_w = a_0$, $b_w = b_0$ and $x_w = x_0$).

In total, we have constructed a triangulation $K$ of $\D$ and a triangulation $T$ of $Q$ with the following properties:
\[
 \tilde{K}^{0} = \{a_0, \ldots, a_{w-1}, b_0, \ldots, b_{w-1} \}, \qquad \tilde{T}^{0} = \{h(a_0), \ldots, h(a_{w-1}), x_0, \ldots, x_{w-1} \},
\]
and the bijection $\varrho : \tilde{K}^0 \to \tilde{T}^0$ defined by $\varrho (a_i) = h(a_i)$ and $\varrho (b_i) = x_i$ for all $i \in \{0, \ldots, w-1\}$ satisfies that for $a,b,c \in \tilde{K}^0$, we have
\[
 abc \in K^2 \quad \mbox{if and only if} \quad \varrho(a) \varrho(b) \varrho (c) \in T^2 .
\]
Therefore, the piecewise affine function $f$ over $K$ such that $f(a) = \varrho (a)$ for all $a \in \tilde{K}^0$ is a homeomorphism from $\D$ onto $Q$ that extends $h$. This finishes the sketch of the construction.

The rest of the section consists in a detailed explicit description of this construction with estimates of all the lengths and angles of the triangulation $K$.
We will prove, thus, the following result.

\begin{theorem}\label{th:Gupta}
Let $\D \subset \R^2$ be a triangle.
Let $\t \in (0,\pi/3]$ be such that
\begin{equation}\label{eq:geqsint}
 \sin \f \geq \sin \t \quad \mbox{for all angles } \ \f \ \mbox{ of the triangle } \ \D.
\end{equation}
Let $m_1$ and $m_2$ be, respectively, a positive lower and an upper bound on the side lengths of $\D$.
Let $0 < l_1 \leq l_2$.
Let $M$ be a $1$-dimensional complex such that
\begin{align}
 & \bigcup M = \p \D , \nonumber \\
 & l_1 \leq |e| \leq l_2 \quad \mbox{for all} \quad e \in M^1. \label{eq:leL}
\end{align}
Let $h: \p\D \to \R^2$ be a homeomorphism that is piecewise affine over $M$.

Then there exist a triangulation $K$ of $\D$, and a homeomorphism $f: \D \to \R^2$ piecewise affine over $K$ such that
\begin{align*}
 & M \subset K , \qquad \tilde{K}^0 \cap \p \D = \tilde{M}^0, \qquad \{ e \in K^1 : e \subset \p \D \} = M^1, \qquad f|_{\p \D} = h \\
 & \sin \f \geq C_0 l_1^{c_{01}} m_1^{c_{02}} m_2^{-c_{04}} (\sin \t)^{c_{03}} \quad \mbox{for all angles } \ \f \ \mbox{ of all triangles in } \ K^2 , \\
 & C_1 l_1^{c_{11}} m_1^{c_{12}} m_2^{-c_{14}} (\sin \t)^{c_{13}} \leq |e| \leq m_2 \quad \mbox{for all } \ e \in K^1 ,
\end{align*}
where
\begin{equation}\label{eq:C0C1}
\begin{split}
 & C_0 := \frac{1}{144}, \quad c_{01} := 2, \quad c_{02} := 2, \quad c_{03} := 4, \quad c_{04} := 4 , \\
 & C_1 := \frac{1}{12}, \quad c_{11} := 1, \quad c_{12} := 1, \quad c_{13} := 2, \quad c_{14} := 1 .
\end{split}
\end{equation}
\end{theorem}

The rest of the section consists of the proof of Theorem \ref{th:Gupta}.

Besides the notation of Theorem \ref{th:Gupta}, in order to do the construction we need to introduce more
notation concerning $\D$.
Let $0 < u < v < w$ be natural numbers.
Let $a_0, a_u, a_v$ be the three vertices of $\D$.
Define $a_w := a_0$.
Let $A,B \in (0, \pi)$ be the angles of $\D$ at the vertices $a_0, a_u$.
Consider points
\begin{eqnarray*}
 & a_0 < a_1 < \cdots < a_{u-1} < a_u & \mbox{ in the order of } \ [a_0, a_u], \\
 & a_u < a_{u+1} < \cdots < a_{v-1} < a_v & \mbox{ in the order of } \ [a_u, a_v], \\
 & a_v < a_{v+1} < \cdots < a_{w-1} < a_w & \mbox{ in the order of } \ [a_v, a_w]
\end{eqnarray*}
such that
\begin{equation*}
 \tilde{M}^0 = \{ a_0,\ldots, a_{w-1} \} , \qquad
 M^1 = \left\{ \ov{a_k a_{k+1}} : k\in \{ 0,\ldots,w-1 \} \right\} .
\end{equation*}
Define $a_{-1} := a_{w-1}$.
Let $r$ and $o$ denote the inradius and incentre, respectively, of $\D$.
Let $h$ be a parameter to be chosen later such that $0 < h < r$. Let $C$ be the circle of centre $o$ and radius $h$.

We will use the following notation about segments, lines and angles.
If $a,b \in \R^2$ are two different points, $\ov{ab}$ denotes the closed segment with endpoints $a,b$. It means the same as $[a,b]$, except that the set $[a,b]$ is equipped with a total order. The length of $\ov{ab}$ is denoted by $|\ov{ab}|$. The straight line passing through $a,b$ is denoted by $\un{ab}$. If $a,b,c \in \R^2$ satisfy $b \notin \{a,c\}$, then $\an{abc}$ denotes the non-oriented angle in $[0,\pi]$ with vertex in $b$ defined by the segments $\ov{ab}$ and $\ov{bc}$. Of course, $\an{abc} = \an{cba}$.

In the construction of this section, it is important to recall the concepts of well-defined polygon and $n$-gon introduced in Section \ref{se:notation}.
In particular, given a natural number $n \geq 3$ and $n$ different points $a_1, \ldots, a_n \in \R^2$, we will say that $a_1 \cdots a_n$ is a \emph{well defined convex $n$-gon} if it is a well defined $n$-gon and is convex as a set.

This section consist of two subsections. Following the notation of Theorem \ref{th:Gupta}, in Subsection \ref{subse:construction} we construct the triangulation $K$ and the homeomorphism $f$, whereas in Subsection \ref{subse:estimates} we estimate the side lengths and angles of $K$.

\subsection{Construction}\label{subse:construction}

This subsection constructs the triangulation $K$ and the homeomorphism $f$ described in Theorem \ref{th:Gupta}.

\begin{figure}
\begin{center}
\input{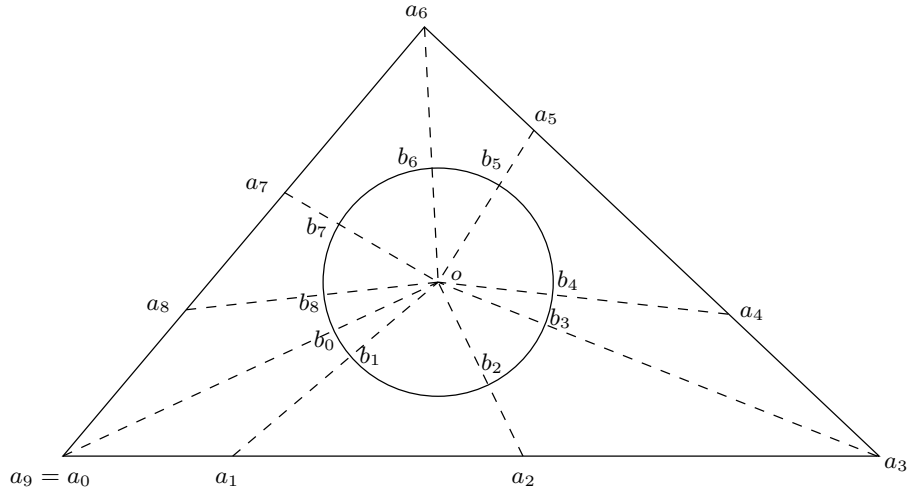}
\caption{\label{fig:Gupta}Points $a_k$ and $b_k$ for $k \in \{0,\ldots,w-1\}$.}
\end{center}
\end{figure}

For $k\in \{ 0, \ldots, w-1 \}$, let $b_k$ be the only point of $\ov{a_k o} \cap C$;
this is represented in Figure \ref{fig:Gupta}, in the particular case when $u=3$, $v=6$ and $w=9$.
Now call $\D' := b_0 \cdots b_{w-1}$.

\begin{lemma}\label{le:D'convex}
$\D'$ is a well-defined convex $w$-gon contained in $\inte{\D}$.
For each $k\in \{ 0, \ldots, w-1 \}$, the quadrilateral $C_k := a_k a_{k+1} b_{k+1} b_k$ is a well-defined convex quadrilateral.
\end{lemma}
\begin{proof}
The geometry of $\D$ implies that the map $\Pi: C \to \p \D$ that sends a point $x\in C$ to the only point of
\[
 \{ o + t (x-o) : t>0 \} \cap \p \D
\]
is a homeomorphism. This implies that $\D'$ is a well-defined polygon. In addition, as any closed disc is strictly convex as a set, then $\D'$ is a well-defined convex $w$-gon. Therefore,
\begin{equation}\label{eq:D'subsetD}
 \D' = \co \{b_0, \ldots, b_{w-1}\} \subset \co C \subset \inte{\D} ;
\end{equation}
here $\co$ denotes the convex hull of a set.

Let $k\in \{ 0, \ldots, w-1 \}$.
To check that $C_k$ is a well-defined quadrilateral that is convex, we show that for every two consecutive vertices $a,b$ of $C_k$, the other two lie on the same connected component of
\[
 \R^2 \setminus \un{ab}.
\]
By \eqref{eq:D'subsetD}, $b_k, b_{k+1} \in \inte{\D}$, and, hence, both  $b_k$ and $b_{k+1}$ lie on the same connected component of
\[
 \R^2 \setminus \un{a_k a_{k+1}}.
\]
Let $H \subset \R^2$ be the open half space such that
\[
 \p H = \un{a_k b_k} \quad \mbox{and} \quad a_{k+1} \in H .
\]
Then $\Pi (C \cap H) = \p \D \cap H$ and, hence, $b_{k+1} \in H$.
Analogously, the two consecutive vertices $a_{k+1}, b_{k+1}$ also enjoy the same property.
Finally, we prove that property for the consecutive vertices $b_k$ and $b_{k+1}$. Since there is a side of $\D$ containing $a_k$ and $a_{k+1}$, then the points $b_k, o, b_{k+1}$ are not aligned. This, together with the facts that $b_k \in (a_k, o)$ and $b_{k+1} \in (a_{k+1}, o)$ imply that $a_k$ and $a_{k+1}$ belong to the same connected component of
\[
 \R^2 \setminus \un{b_k b_{k+1}} .
\]
This proves that $C_k$ is a well-defined convex quadrilateral.
\end{proof}

\begin{figure}
\begin{center}
\input{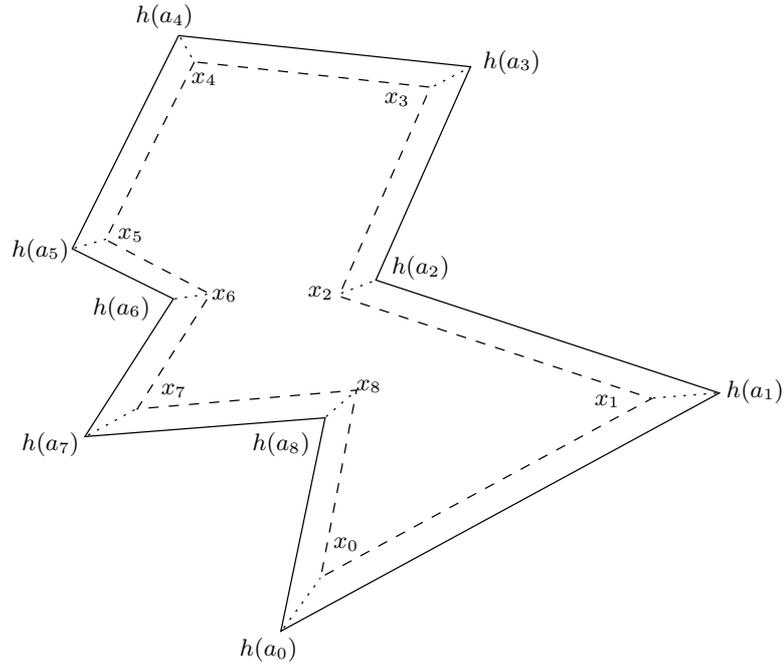}
\caption{\label{fig:x0x9}Polygons $Q$ and $Q'$.}
\end{center}
\end{figure}

As $h$ is a homeomorphism piecewise affine over $M$, then $h(\p \D)$ equals the boundary of the polygon $Q := h(a_0) \cdots h(a_{w-1})$. The construction of Lemma \ref{le:x0x8} below is represented in Figure \ref{fig:x0x9}, in the particular case when $u=3$, $v=6$ and $w=9$.

\begin{lemma}\label{le:x0x8}
There exist $x_0, \ldots, x_{w-1} \in Q$ such that
$Q' := x_0 \cdots x_{w-1}$ is a well-defined $w$-gon contained in $\inte{Q}$ and that for each $k\in \{ 0, \ldots, w-1 \}$ the quadrilateral $Q_k := h(a_k) h(a_{k+1}) x_{k+1} x_k$ is a well-defined convex quadrilateral, contained in $Q \setminus \inte{Q'}$, and $Q' \cap \inte{Q_k} = \varnothing$.
Moreover, for each $0 \leq j < k \leq w-1$, the set $Q_j \cap Q_k$ is either empty or a common side of $Q_j$ and $Q_k$.
Finally, 
\[
 Q = Q' \cup \bigcup_{k=0}^{w-1} Q_k .
\]
\end{lemma}
\begin{proof}
For every $k\in \{ 0, \ldots, w-1 \}$, let $e_k \in \R^2$ be the unit vector such that
\begin{equation}\label{eq:bisects}
 \un{h(a_k) [h(a_k) + e_k]} \quad \mbox{bisects} \quad \an{\scriptstyle h(a_{k-1}) h(a_k) h(a_{k+1})}
\end{equation}
and such that $h(a_k) + t e_k \in \inte{Q}$ for every $t>0$ sufficiently small.

Fix $r>0$.
For each $k\in \{ 0, \ldots, w-1 \}$, define $t_k >0$ and $x_k$ through the property that
\[
 x_k := h(a_k) + t_k e_k
\]
satisfy $\dist(x_k, \un{h(a_k) h(a_{k+1})}) = r$. Property \eqref{eq:bisects} implies that in fact
\[
 \dist(x_k, \un{h(a_k) h(a_{k+1})}) =  \dist(x_k, \un{h(a_{k-1}) h(a_k)}) = r
\]
and, hence, the lines
\[
 \un{h(a_k) h(a_{k+1})} \quad \mbox{and} \quad \un{x_k x_{k+1}} \quad \mbox{are parallel.}
\]
If $r$ is small then, for each $k\in \{ 0, \ldots, w-1 \}$,
\[
 x_k \in \inte{Q} \quad \mbox{and} \quad \ov{h(a_k) x_k} \cap \ov{h(a_{k+1}) x_{k+1}} = \varnothing ;
\]
this and the choice of $e_k$ imply that $Q_k$ is a well-defined
quadrilateral, in fact, a trapezoid and hence convex.
Moreover, if $r$ is small then $Q_k \subset Q \setminus \inte{Q'}$,
and $Q'$ is a well-defined polygon contained in $\inte{Q}$.
Furthermore, if $r$ is small then for each $0 \leq j < k \leq w-1$,
the set $Q_j \cap Q_k$ is either empty or a common side of $Q_j$ and
$Q_k$, and $Q' \cap \inte{Q_k} = \varnothing$. If $Q'$ is a $w$-gon
then we are done. If not, we can perturb the vertices $x_0, \ldots,
x_{w-1}$ so that all properties of the statement are true, including
the fact that $Q'$ is a $w$-gon.

Finally, the equality $Q = Q' \cup \bigcup_{k=0}^{w-1} Q_k$ is obvious from the construction.
\end{proof}

Choose any $x_0, \ldots, x_{w-1} \in \inte{Q}$ satisfying the conditions of the statement of Lemma \ref{le:x0x8}.
It is a well-known result (e.g., the proof of Lemma 3.3.2 in Dettman \cite{Dettman}) that every polygon can be triangulated without adding extra vertices; therefore,
there exists a triangulation $N$ of $Q'$ such that $\tilde{N}^0 = \{ x_0, \ldots, x_{w-1} \}$.
Define $\varrho : \{ b_0, \ldots, b_{w-1} \} \to \tilde{N}^0$ through $\varrho (b_k) = x_k$ for all $k\in \{ 0, \ldots, w-1 \}$.
Now define the $2$-dimensional complex $L$ as follows:
\begin{align*}
 & \tilde{L}^0 := \varrho^{-1} (\tilde{N}^0), \\
 & L^1 := \left\{ \ov{ab} : \quad a,b \in \tilde{L}^0, \quad \ov{\varrho(a) \varrho(b)} \in N^1 \right\}, \\
 & L^2 := \left\{ abc : \quad a,b,c \in \tilde{L}^0, \quad \varrho(a) \varrho(b) \varrho(c) \in N^2 \right\} .
\end{align*}
The following lemma is based on an observation by Aronov, Seidel and Souvaine \cite{ASR},
according to which every triangulation of an $n$-gon that does not add any extra vertex induces a triangulation of any well-defined convex $n$-gon.
\begin{lemma}\label{le:L}
$L$ is a triangulation of $\D'$ such that $\ov{b_k b_{k+1}} \in L^1$ for all $k \in \{ 0, \ldots, w-1 \}$.
\end{lemma}
\begin{proof}
Let $\s \in L^2$.
Then $\s = abc$ for three different points $a,b,c \in \tilde{L}^0$.
Since, by Lemma \ref{le:D'convex}, $\D'$ is a well-defined convex $w$-gon, then $\s$ is a triangle contained in $\D'$. In particular, $\bigcup L^2 \subset \D'$.

We prove that the sides of any $\s \in L^2$ belong to $L^1$. Indeed, let $\s \in L^2$. Then $\s = abc$ for some $a,b,c \in \tilde{L}^0$ with $\varrho(a) \varrho(b) \varrho(c) \in N^2$. Since $N$ is a complex then $\ov{\varrho(a) \varrho(b)} \in N^1$, and, hence, $\ov{ab} \in L^1$. Analogously, $\ov{bc}, \ov{ca}\in L^1$.
Similarly, one proves that the endpoints of any element of $L^1$ belong to $\tilde{L}^0$.

The last paragraph shows in particular that
\[
 \bigcup L^2 = \bigcup L = \bigcup_{\s \in L^2} \inte{\s} \cup \bigcup_{e \in L^1} \inte{e} \cup \tilde{L}^0 .
\]
We are going to prove that $\bigcup L$ is open in $\D'$;
to do this, we show that every element of 
\[
  \inte{\s} \quad \mbox{(for } \ \s \in L^2 \mbox{ )} , \qquad 
 \inte{e} \quad \mbox{(for } \ e \in L^1 \mbox{ )} , \qquad
 \tilde{L}^0
\]
has a neighborhood $N$ such that $N \cap \D' \subset \bigcup L$.
If $x \in \inte{\s}$ for some $\s \in L^2$ then there exists a neighbourhood $N$ of $x$ such that $N \subset \s \subset \bigcup L$.
If $x \in \inte{e}$ for some $e \in L^1$ then there exist $j\neq k \in \{0,\ldots,w-1\}$ such that $e = \ov{b_j b_k}$ and $\ov{x_j x_k} \in N^1$. Suppose $e$ is a side of $\D'$; then $j$ is congruent to $k-1$ or $k+1$ modulo $w$, and $\ov{x_j x_k}$ is a side of $Q'$. As $N$ is a triangulation of $Q'$ such that $\tilde{N}^0 = \{ x_0, \ldots, x_{w-1} \}$, there exists $i \in \{0, \ldots,w-1\}$ such that $x_j x_k x_i \in N^2$. Then $b_j b_k b_i \in L^2$ and there exists a neighbourhood $N$ of $x$ such that $N \cap \D'\subset b_j b_k b_i \subset \bigcup L$. Suppose $e$ is not a side of $\D'$; then $j$ is not congruent to $k-1$ or $k+1$ modulo $w$, and $\ov{x_j x_k}$ is not a side of $Q'$. As $N$ is a triangulation of $Q'$ such that $\tilde{N}^0 = \{ x_0, \ldots, x_{w-1} \}$, there exist $i\neq i' \in \{0, \ldots,w-1\}$ such that $x_j x_k x_i, x_j x_k x_{i'} \in N^2$. Then $b_j b_k b_i, b_j b_k b_{i'} \in L^2$ and there exists a neighbourhood $N$ of $x$ such that $N \subset b_j b_k b_i \cup b_j b_k b_{i'} \subset \bigcup L$.
Finally, if $x \in \tilde{L}^0$ define $S_N$ as the set of $\s \in N^2$ such that $\varrho (x) \in \s$, and $S_L$ as the set of $\s \in L^2$ such that $\s = abc$ for some $a,b,c \in \tilde{L}^0$ with $\varrho(a) \varrho(b) \varrho(c) \in S_N$. Then there exists a neighbourhood $N$ of $x$ such that $N \cap \D' \subset \bigcup S_L \subset \bigcup L$.
This proves that $\bigcup L$ is open in $\D'$. Clearly, $\bigcup L$ is closed as a finite union of closed sets. Therefore, $\bigcup L$ is non-empty, open and closed in the connected set $\D'$, whence $\bigcup L = \D'$.

Now let $e \neq c \in L^1$ satisfy $e \cap c \neq \varnothing$. Without loss of generality, suppose that $e = \ov{b_0 b_i}$ and $c = \ov{b_j b_k}$ for some $i,j,k \in \{1,\ldots,w-1\}$ with $j <k$. As $e \cap c \neq \varnothing$ and $\D'$ is a well-defined convex $w$-gon, then $j \leq i$. If $j=i$ then $e \cap c \in L^0$. If $j < i$ then $\ov{x_0 x_i} \cap \ov{x_j x_k} \in \inte{Q'}$ with $\ov{x_0 x_i}, \ov{x_j x_k} \in N^1$, which contradicts the fact that $N$ is a triangulation of $Q'$ such that $\tilde{N}^0 = \{ x_0, \ldots, x_{w-1} \}$. Therefore, $e \cap c \in L^0$.
Since $\tilde{L}^0 = \{ b_0, \ldots, b_{w-1}\}$, this is enough to prove that $L$ is a complex.

Finally, since $N$ is a triangulation of $Q'$ such that $\tilde{N}^0 = \{ x_0, \ldots, x_{w-1} \}$, for each $k \in \{ 0, \ldots, w-1 \}$ there exists a triangle in $N^2$ one of which sides is $\ov{x_k x_{k+1}}$. Since $N$ is a triangulation, necessarily $\ov{x_k x_{k+1}} \in N^1$. Therefore, $\ov{b_k b_{k+1}} \in L^1$.
\end{proof}

\begin{figure}
\begin{center}
\includegraphics{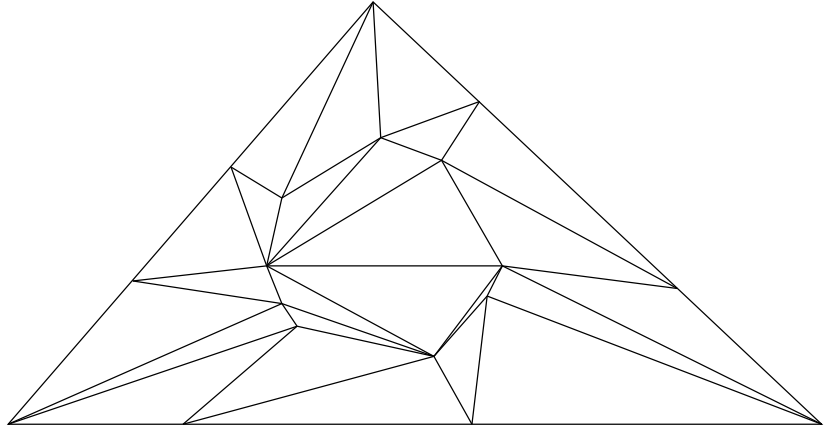}
\caption{\label{fig:Gupta2}Construction of the triangulation $K$.}
\end{center}
\end{figure}

Define the $2$-dimensional complex $K$ as follows:
\begin{align*}
 & K^0 := M^0 \cup L^0 , \\
 & K^1 := L^1 \cup \bigcup_{k=0}^{w-1} \{ \ov{a_k a_{k+1}} , \ov{b_k b_{k+1}} , \ov{a_k b_k} , \ov{a_k b_{k+1}} \} , \qquad
  K^2 := L^2 \cup \bigcup_{k=0}^{w-1} \{ a_k a_{k+1} b_{k+1} , a_k b_{k+1} b_k \} .
\end{align*}
This complex $K$ is represented in Figure \ref{fig:Gupta2}, in the particular case when $u=3$, $v=6$ and $w=9$.

\begin{lemma}\label{le:K}
$K$ is a triangulation of $\D$.
\end{lemma}
\begin{proof}
By Lemmas \ref{le:L} and \ref{le:D'convex}, $\bigcup L = \D' \subset \D$.
Moreover, by construction and Lemma \ref{le:D'convex}, $a_k, b_k \in \D \cup \D' = \D$ for each $k \in \{ 0, \ldots, w-1 \}$. As $\D$ is convex then $\bigcup K^2 \subset \D$. Now we note that, by Lemmas \ref{le:L} and \ref{le:D'convex},
\[
 \bigcup K = \D' \cup \bigcup_{k=0}^{w-1} C_k =
 \inte{\D'} \cup \bigcup_{k=0}^{w-1} \left[ (a_k, a_{k+1}) \cup (b_k, b_{k+1}) \cup (a_k, b_k) \cup \{ a_k, b_k \} \right] .
\]
We are going to prove that $\bigcup K$ is open in $\D$;
to do this, we show that every element of 
\[
 \inte{\D'} , \quad (a_k, a_{k+1}) , \quad (b_k, b_{k+1}) , \quad (a_k, b_k) , \quad \{ a_k, b_k \} ,
\]
(for $k \in \{0,\ldots,w-1\}$) has a neighbourhood $N$ such that $N \cap \D \subset \bigcup K$.
Indeed, the following properties are easy to verify.
If $x \in \inte{\D'}$ then there exists a neighbourhood $N$ of $x$ such that $N \subset \D' \subset \bigcup K$.
Let $k \in \{ 0, \ldots, w-1 \}$.
If $x \in (a_k, a_{k+1})$ then there exists a neighbourhood $N$ of $x$ such that $N \cap \D \subset a_k a_{k+1} b_{k+1} \subset \bigcup K$.
If $x \in (b_k, b_{k+1})$ then there exists a neighbourhood $N$ of $x$ such that $N  \subset \D' \cup a_k b_{k+1} b_k \subset \bigcup K$.
If $x \in (a_k, b_k)$ then there exists a neighbourhood $N$ of $x$ such that $N \subset a_{k-1} a_k b_k \cup a_k b_{k+1} b_k \subset \bigcup K$.
There exists a neighbourhood $N$ of $a_k$ such that $N \cap \D \subset a_{k-1} a_k b_k \cup C_k \subset \bigcup K$.
There exists a neighbourhood $N$ of $b_k$ such that $N \subset C_{k-1} \cup a_k b_{k+1} b_k \cup \D' \subset \bigcup K$.
This proves that $\bigcup K$ is open in $\D$. Clearly, $\bigcup K$ is closed as a finite union of closed sets. Therefore, the non-empty set $\bigcup K$ is open and closed in the connected set $\D$, which implies $\bigcup K = \D$.

All elements of $K^2$ are triangles. Indeed, by Lemma \ref{le:L}, every element of $L^2$ is a triangle. In addition, by Lemma \ref{le:D'convex},
for each $k\in \{ 0, \ldots, w-1 \}$, the quadrilateral $C_k$ is a well-defined convex quadrilateral, and, hence, $a_k a_{k+1} b_{k+1}$ and $a_k b_{k+1} b_k$ are triangles whose union is $C_k$.

By Lemma \ref{le:L} and the definition of $K$, the sides of any $\s \in K^2$ belong to $K^1$, and
the endpoints of any $e \in K^1$ belong to $\tilde{K}^0$.

Finally, let $\s \neq \tau \in K^2$ satisfy $\s \cap \tau \neq \varnothing$; we want to prove that $\s \cap \tau \in K^0 \cup K^1$.
If $\s, \tau \in L^2$ then, by Lemma \ref{le:L}, $\s \cap \tau \in L^0 \cup L^1 \subset K^0 \cup K^1$.
Let $j,k \in \{ 0, \ldots, w-1 \}$.
If $\s \in L^2$ and $\tau = a_k a_{k+1} b_{k+1}$ then $\s \cap \tau = \{ b_{k+1} \} \in L^0 \subset K^0$.
If $\s \in L^2$ and $\tau = a_k b_{k+1} b_k$ then, by Lemma \ref{le:L}, $\s \cap \tau \in \{ \{b_k\} , \{b_{k+1}\} , \ov{b_k b_{k+1}} \} \subset L^0 \cup L^1 \subset K^0 \cup K^1$.
If $\s = a_k a_{k+1} b_{k+1}$ and $\tau = a_j a_{j+1} b_{j+1}$ then $\s \cap \tau \in \{ \{a_k\} , \{a_{k+1}\} \} \subset M^0 \subset K^0$.
If $\s = a_k a_{k+1} b_{k+1}$ and $\tau = a_j b_{j+1} b_j$ then $\s \cap \tau \in \{ \ov{a_k b_{k+1}} , \ov{a_{k+1} b_{k+1}} \} \subset K^1$.
If $\s = a_k b_{k+1} b_k$ and $\tau = a_j b_{j+1} b_j$ then $\s \cap \tau \in \{ \{b_k\} , \{b_{k+1}\} \} \subset L^0 \subset K^0$. This completes the proof.
\end{proof}

By Lemma \ref{le:K}, there exists a unique continuous function $f: \D \to \R^2$ piecewise affine over $K$ such that
\begin{equation*}
 f(a_k) = h(a_k) , \quad f(b_k) = x_k , \qquad k\in \{ 0, \ldots, w-1 \} .
\end{equation*}
\begin{lemma}
$f$ is a homeomorphism onto $Q$ that coincides with $h$ in $\p \D$.
\end{lemma}
\begin{proof}
Since $\D$ is compact and $f$ is continuous, to prove that $f$ is a homeomorphism it suffices to show that $f$ is injective. By Lemma \ref{le:x0x8} and the fact that $h$ is a homeomorphism,
\begin{equation}\label{eq:finjectiveK0}
 f \ \mbox{ is injective in  } \ \tilde{K}^0.
\end{equation}
Now we prove that $f$ is injective in $\bigcup K^1$. For each $e \in
K^1$ the function $f|_e$ is affine and (by \eqref{eq:finjectiveK0})
non constant, hence injective; in particular, $f(\inte{e})$ equals
the interior of $f(e)$. Let $e\neq c \in K^1$. If $e,c \in L^1$ then
$f(e), f(c) \in N^1$ and (by \eqref{eq:finjectiveK0}) $f(e) \neq
f(c)$; since $N$ is a complex then $f(\inte{e}) \cap f(c) =
\varnothing$. If $e \in L^1$ and $c \notin L^1$ then by Lemma
\ref{le:x0x8}, $f(\inte{e}) \subset \inte{Q'}$ and $f(c) \cap
\inte{Q'} = \varnothing$; therefore, $f(\inte{e}) \cap f(c) =
\varnothing$; similarly, one proves $f(e) \cap f(\inte{c}) =
\varnothing$. If $e,c \notin L^1$ then, by Lemma \ref{le:x0x8},
\[
 f(e) \cap f(c) \in N^0 \cup \left\{ \{h(a_0)\}, \ldots, \{h(a_{w-1})\} , \varnothing \right\}
\]
and, hence, $f(\inte{e}) \cap f(c) =
\varnothing$. As \eqref{eq:finjectiveK0}, this proves that $f$ is
injective in $\bigcup K^1$.

Now we prove that $f$ is injective in each $\s \in K^2$.
If $\s \in L^2$ then $\s = abc$ for some $a,b,c \in \tilde{L}^0$ such that $f(\s) = \varrho(a) \varrho(b) \varrho(c) \in N^2$; as $N$ is a triangulation then $f(\s)$ is a triangle and, hence, $f|_{\s}$ is injective.
If $\s \notin L^2$ then, by Lemma \ref{le:x0x8}, $f(\s)$ is a triangle and, hence, $f|_{\s}$ is injective.

Now let $\s \neq \tau \in K^2$. To complete the proof of the injectivity of $f$ we only have to show $f(\inte{\s}) \cap f(\tau) = \varnothing$.
If $\s, \tau \in L^2$ then $f(\s)\neq f(\tau) \in N^2$; as $N$ is a triangulation then $f(\inte{\s}) \cap f(\tau) = \varnothing$.
If $\s \in L^2$ and $\tau \notin L^2$ then $f(\s) \in N^2$, hence $f(\s) \subset Q'$, and $f(\tau) \subset Q_k$ for some $k \in \{ 0, \ldots, w-1 \}$; by Lemma \ref{le:x0x8}, $f(\inte{\s}) \cap f(\tau) = f(\s) \cap f(\inte{\tau}) = \varnothing$.
If $\s, \tau \notin L^2$ then, by Lemma \ref{le:x0x8}, $f(\inte{\s}) \cap f(\tau) = \varnothing$.
This completes the proof that $f$ is a homeomorphism.

Now, by Lemmas \ref{le:K} and \ref{le:x0x8},
\begin{align*}
 f(\D)
 & = \bigcup_{\s \in K^2} f(\s)
 = \bigcup_{\s \in L^2} f(\s) \cup \bigcup_{k=0}^{w-1} f(a_k a_{k+1} b_{k+1}) \cup f(a_k b_{k+1} b_k) \\
 & = \bigcup N^2 \cup \bigcup_{k=0}^{w-1} h(a_k) h(a_{k+1}) x_{k+1} \cup h(a_k) x_{k+1} x_k
 = Q' \cup \bigcup_{k=0}^{w-1} Q_k = Q .
\end{align*}

Finally, we note that
\[
 \p\D = \bigcup_{k=0}^{w-1} \ov{a_k a_{k+1}} = \bigcup M
\]
and that both $f|_{\p\D}$ and $h$ are piecewise affine over $M$ and coincide in $\tilde{M}^0$. Therefore,
$f$ coincides with $h$ in $\p \D$.
\end{proof}

\subsection{Estimates}\label{subse:estimates}

Our aim in this subsection is to give upper and lower bounds of
\[
 |e| \qquad \mbox{for} \quad e \in K^1 ,
\]
and lower bounds of
\[
 \sin \f \qquad \mbox{for all angles } \ \f \ \mbox{ of all triangles in } \ K^2
\]
in terms of $l_1,l_2,\t, m_1, m_2$.
Of course, $K$ is the complex constructed in Subsection \ref{subse:construction},
and we are following the notation of Theorem \ref{th:Gupta}.
In Subsection \ref{subse:construction} we took any $h$ in $(0,r)$, where $r$ is the inradius of the triangle $\D$.
In this subsection we make the choice
\[
 h := \frac{1}{6} m_1 \sin \t.
\]
The fact that this $h$ is indeed less than $r$ is proved in Lemma \ref{le:inradiusineq} below.
With this choice, we shall estimate the following lengths:
\begin{align*}
 & |\ov{a_k a_{k+1}}| , \ |\ov{a_k b_k}| , \ |\ov{a_k b_{k+1}}| ,
 \qquad k \in \{ 0, \ldots, w-1 \} , \\
 & |\ov{b_i b_j}| , \qquad i,j \in \{ 0, \ldots, w-1 \}, \quad i < j ,
\end{align*}
and the sine of the following angles:
\begin{align*}
 & \an{a_k a_{k+1} b_{k+1}} , \ \an{a_{k+1} b_{k+1} a_k} , \ \an{b_{k+1} a_k a_{k+1}} , \
 \an{a_k b_{k+1} b_k} , \ \an{b_{k+1} b_k a_k} , \ \an{b_k a_k b_{k+1}} ,
 \quad k \in \{ 0, \ldots, w-1 \} , \\
 & \an{b_i b_j b_k} , \qquad i,j,k \in \{ 0, \ldots, w-1 \}, \quad i < k, \quad j \notin \{i,k\} .
\end{align*}
Before doing that, we need some trigonometric inequalities and identities. The following two are elementary:
\begin{align}
 & \frac{\sin x}{2} \leq \sin \frac{x}{2} , \qquad x \in [0, 2 \pi] , \label{eq:ineqsin}\\
 & 2 \sin ^2 x = 1 - \cos 2x, \qquad x \in \R . \label{eq:identitysin}
\end{align}
In the following lemma we estimate the inradius of a triangle from below.
\begin{lemma}\label{le:inradiusineq}
Consider a triangle $\D$ with inradius $r$, minimum side length $m_1$, and $\t \in (0 , \pi/3]$ satisfies \eqref{eq:geqsint}. Then
\[
 r \geq \frac{1}{3} m_1 \sin \t .
\]
\end{lemma}
\begin{proof}
Let $|\D|$ denote the area of $\D$, and $s$ its semiperimeter. Then, by elementary trigonometry,
$r = |\D|/s$.
Now let $a\leq b\leq c$ be the three side lengths of the triangle. By elementary geometry and \eqref{eq:geqsint}, we have $2 |\D| \geq b c \sin \t$ and $2 s \leq 2b + c$.
Therefore,
\[
 r \geq \frac{b}{3} \sin \t \geq \frac{m_1}{3} \sin \t ,
\]
as required.
\end{proof}

In the following Lemmas we estimate the lengths and angles of the triangulation $K$ constructed in Subsection \ref{subse:construction}.

\begin{lemma}\label{le:akak+1}
For all $k \in \{ 0, \ldots, w-1 \}$,
\[
 l_1 \leq |\ov{a_k a_{k+1}}| \leq l_2 .
\]
\end{lemma}
\begin{proof}
This is immediate due to inequalities \eqref{eq:leL}.
\end{proof}

\begin{lemma}\label{le:akbk}
For all $k \in \{ 0, \ldots, w-1 \}$,
\[
 \frac{1}{6} m_1 \sin \t \leq |\ov{a_k b_k}| \leq m_2 .
\]
\end{lemma}
\begin{proof}
Fix $k \in \{ 0, \ldots, w-1 \}$.
As $b_k \in (a_k, o)\cap C$, by the geometry of $\D$,
\[
 |\ov{a_k b_k}| = |\ov{a_k o}| - h \geq r - h .
\]
Lemma \ref{le:inradiusineq} concludes the first inequality of the
statement. The second inequality is obvious, since $|\ov{a_k b_k}|
\leq \diam \D \leq m_2$.
\end{proof}

\begin{lemma}\label{le:bibj}
For all $0 \leq i < j \leq w-1$,
\[
 \frac{l_1 m_1 \sin^2 \t}{12 m_2} \leq |\ov{b_i b_j}| \leq \frac{1}{3} m_1 \sin \t .
\]
\end{lemma}
\begin{proof}
For $0 \leq i < j \leq w-1$ we have $b_i, b_j \in C$ and, hence, $|\ov{b_i b_j}| \leq \diam C = 2h$. This proves the second inequality of the statement. Now,
by the Law of Cosines applied to the triangle $b_i o b_j$,
\begin{equation*}
 |\ov{b_i b_j}|^2 = 2 h^2 \left( 1 - \cos \an{b_i o b_j} \right) .
\end{equation*}
This equality and the following one
\begin{equation}\label{eq:minbiobj}
  \min_{0 \leq i < j \leq w-1} \an{b_i o b_j} = \min_{0 \leq k \leq w-1} \an{b_k o b_{k+1}}
\end{equation}
imply
\begin{equation}\label{eq:minbibj}
 \min_{0 \leq i < j \leq w-1} |\ov{b_i b_j}|^2 = \min_{0 \leq k \leq w-1} |\ov{b_k b_{k+1}}|^2 = 2 h^2 \min_{0 \leq k \leq w-1} \left( 1 - \cos \an{b_k o b_{k+1}} \right) .
\end{equation}

Now let $0 \leq k \leq u-1$. As $b_k \in (a_k , o)$ and $b_{k+1} \in (a_{k+1} , o)$ then
$\an{b_k o b_{k+1}} = \an{a_k o a_{k+1}}$.
By the Law of Sines applied to $a_k o a_{k+1}$,
\[
 \sin \an{a_k o a_{k+1}} = \frac{|\ov{a_k a_{k+1}}|}{|\ov{a_{k+1} o}|} \sin \an{o a_k a_{k+1}},
\]
but, thanks to the geometry of $\D$, inequalities \eqref{eq:geqsint} and \eqref{eq:ineqsin}, and the facts that $\un{a_0 o}$ bisects $\an{a_v a_0 a_u}$, and $\un{a_u o}$ bisects $\an{a_0 a_u a_v}$, we have
\[
 \sin \an{o a_k a_{k+1}} \geq \min \{ \sin \an{o a_0 a_u} , \sin \an{o a_u a_0} \} =
 \min \{ \sin \frac{A}{2} , \sin \frac{B}{2} \} \geq \sin \frac{\t}{2} \geq \frac{\sin \t}{2} .
\]
Since we have $|\ov{a_k a_{k+1}}| \geq l_1$ and $|\ov{a_{k+1} o}| \leq m_2$, then, also by \eqref{eq:ineqsin},
\begin{equation}\label{eq:akoak+12}
  \sin \an{a_k o a_{k+1}} \geq \frac{l_1 \sin \t}{2 m_2} \quad \mbox{and} \quad
 \sin \frac{\an{a_k o a_{k+1}}}{2} \geq \frac{l_1 \sin \t}{4 m_2} ;
\end{equation}
this, together with \eqref{eq:identitysin}, \eqref{eq:minbibj} and the symmetry of the argument, completes the proof.
\end{proof}

\begin{lemma}\label{eq:bibjbk}
For all $i,j,k \in \{ 0, \ldots, w-1 \}$ with $i < k$ and $j \notin \{i,k\}$,
\[
 \sin \an{b_i b_j b_k} \geq \frac{l_1 \sin \t}{4 m_2} .
\]
\end{lemma}
\begin{proof}
Let $k \in \{ 0, \ldots, w-1 \}$.
By elementary plane geometry, for any $x \in C \setminus \{ b_k, b_{k+1} \}$,
\[
 \an{b_k x b_{k+1}} = \frac{1}{2} \an{b_k o b_{k+1}}.
\]
As $b_k \in (a_k, o)$ and $b_{k+1} \in (a_{k+1}, o)$ then $\an{b_k o b_{k+1}} = \an{a_k o a_{k+1}}$. Equations \eqref{eq:minbiobj} and \eqref{eq:akoak+12} conclude the proof.
\end{proof}

\begin{lemma}\label{le:akbk+1}
For all $k \in \{ 0, \ldots, w-1 \}$,
\[
 \frac{1}{6} m_1 \sin \t \leq |\ov{a_k b_{k+1}}| \leq m_2 .
\]
\end{lemma}
\begin{proof}
Let $k \in \{ 0, \ldots, w-1 \}$.
Since $b_{k+1} \in C$ and $b_k \in C \cap (a_k,o)$, by the geometry of $\D$ and Lemma \ref{le:akbk},
\[
 |\ov{a_k b_{k+1}}| \geq \min_{x \in C} |\ov{a_k x}| = |\ov{a_k b_k}| \geq \frac{1}{6} m_1 \sin \t .
\]
This shows
the first inequality of the statement. The second inequality is obvious, since $\diam \D \leq m_2$.
\end{proof}

\begin{lemma}\label{le:akak1bk1}
For all $k \in \{ 0, \ldots, w-1 \}$,
\begin{align}
 \sin \an{a_k a_{k+1} b_{k+1}} & \geq \frac{l_1 m_1 \sin^2 \t}{6 l_2 m_2} , \label{eq:akak+1bk+1}\\
 \sin \an{a_{k+1} b_{k+1} a_k} & \geq \frac{l_1^2 m_1 \sin^2 \t}{6 l_2 m_2^2} , \nonumber \\
 \sin \an{b_{k+1} a_k a_{k+1}} & \geq \frac{l_1 m_1^2 \sin^3 \t}{36 l_2 m_2^2} . \nonumber
\end{align}
\end{lemma}
\begin{proof}
Fix $k \in \{ 0, \ldots, w-1 \}$. By the Law of Sines applied to $a_k a_{k+1} o$,
\[
 \sin \an{a_k a_{k+1} o} = \frac{|\ov{a_k o}|}{|\ov{a_k a_{k+1}}|} \sin \an{a_k o a_{k+1}} .
\]
The estimates of \eqref{eq:akoak+12}, Lemmas \ref{le:akak+1} and \ref{le:akbk}, and the facts $b_k \in (a_k, o) \cap C$ and $b_{k+1} \in (a_{k+1}, o)$ demonstrate \eqref{eq:akak+1bk+1}.

By the Law of Sines applied to $a_k a_{k+1} b_{k+1}$,
\[
 \sin \an{a_{k+1} b_{k+1} a_k} = \frac{|\ov{a_k a_{k+1}}|}{|\ov{a_k b_{k+1}}|} \sin \an{a_k a_{k+1} b_{k+1}},  \qquad
 \sin \an{b_{k+1} a_k a_{k+1}} = \frac{|\ov{a_{k+1} b_{k+1}}|}{|\ov{a_k b_{k+1}}|} \sin \an{a_k a_{k+1} b_{k+1}} .
\]
Inequality \eqref{eq:akak+1bk+1} and Lemmas \ref{le:akak+1}, \ref{le:akbk+1} and \ref{le:akbk} show the second and third inequalities of the statement.
\end{proof}

\begin{lemma}\label{eq:akbk1bk}
For all $k \in \{ 0, \ldots, w-1 \}$,
\begin{align}
 \sin \an{a_k b_{k+1} b_k} & \geq \frac{l_1 m_1 \sin^2 \t}{12 m_2^2} , \label{eq:akbk+1bk} \\
 \sin \an{b_{k+1} b_k a_k} & \geq \frac{l_1 m_1^2 \sin^3 \t}{72 m_2^3} , \nonumber \\
 \sin \an{b_k a_k b_{k+1}} & \geq \frac{l_1^2 m_1^2 \sin^4 \t}{144 m_2^4} . \nonumber
\end{align}
\end{lemma}
\begin{proof}
Fix $k \in \{ 0, \ldots, w-1 \}$. By the Law of Sines applied to $a_k b_{k+1} o$,
\[
 \sin \an{o a_k b_{k+1}} = \frac{|\ov{b_{k+1} o}|}{|\ov{a_k b_{k+1}}|} \sin \an{a_k o b_{k+1}} .
\]
The estimates of \eqref{eq:akoak+12} and Lemma \ref{le:akbk+1}, and the facts $b_k \in (a_k, o)$ and $b_{k+1} \in (a_{k+1}, o)\cap C$ demonstrate \eqref{eq:akbk+1bk}.

By the Law of Sines applied to $a_k b_{k+1} b_k$,
\[
 \sin \an{b_{k+1} b_k a_k} = \frac{|\ov{a_k b_{k+1}}|}{|\ov{a_k b_k}|} \sin \an{a_k b_{k+1} b_k} , \qquad
 \sin \an{b_k a_k b_{k+1}} = \frac{|\ov{b_k b_{k+1}}|}{|\ov{a_k b_k}|} \sin \an{a_k b_{k+1} b_k} .
\]
Inequality \eqref{eq:akbk+1bk} and Lemmas \ref{le:akbk}, \ref{le:akbk+1} and \ref{le:bibj} demonstrate the second and third inequalities of the statement.
\end{proof}

From Lemmas \ref{le:akak+1}, \ref{le:akbk}, \ref{le:bibj} and \ref{le:akbk+1} we see that $m_2$ is an upper bound for the lengths of the segments in $K^1$, and
\[
 \frac{l_1 m_1 \sin^2 \t}{12 m_2}
\]
is a lower bound for the lengths of the segments in $K^1$.
From Lemmas \ref{eq:bibjbk}, \ref{le:akak1bk1} and \ref{eq:akbk1bk} we see that
\[
 \frac{l_1^2 m_1^2 \sin^4 \t}{144 m_2^4}
\]
is a lower bound for the sine of the angles of the triangles in $K^2$.
This finishes the proof of Theorem \ref{th:Gupta}.

\section{Construction of the triangulation}\label{se:triangulation}

In this section we construct the triangulation over which our approximating homeomorphism $f$ will be piecewise affine.
The idea is to start with a fine triangulation whose triangles are `almost equilateral and with about the same size', in the sense of Proposition \ref{prop:triangulation} below.
Then we refine the skeleton of that triangulation using Theorem \ref{th:skeleton}; we thus obtain a homeomorphism $g$ piecewise affine over the skeleton of a suitable triangulation that approximates our original homeomorphism $h$.
Finally, we extend $g$ from the boundary of each triangle of the triangulation to the whole triangle, using Theorem \ref{th:Gupta}.
The result will be a piecewise affine homeomorphism $f$ that approximates $h$ in the supremum norm.

We will need the following result, a proof of which can be found in Shewchuk \cite{Shewchuk}. It says that every polygon admits triangulations as fine as we wish whose triangles are `almost equilateral and with about the same size'; in the finite element literature, this is called a \emph{quasiuniform} triangulation.

\begin{proposition}\label{prop:triangulation}
Let $\O \subset \R^2$ be a closed polygon. Then there exist $\e_0 >0$,
\begin{equation}\label{eq:td}
 \t \in ( 0, \pi/3] , \qquad d \in (0,1)
\end{equation}
such that for all $0< \e \leq \e_0$ there exists a triangulation $L$ of $\O$ satisfying
\begin{equation}\label{eq:singeq}
 \sin \f \geq \sin \t \quad \mbox{for all angles } \f \ \mbox{ of all triangles in } \ L^2,
\end{equation}
and
\[
 d \e \leq |e| \leq \e , \qquad e \in L^1.
\]
\end{proposition}

The following property is implicit in Moise \cite{Moise} (proof of Theorem 3, Chapter 6); a proof of a more general result can be found in Paul \cite{Paul} (Theorem 3.3 of Part I).
\begin{lemma}\label{le:fsubsetBh}
Let $\s \subset \R^2$ be a closed triangle. Let $f,h : \s \to \R^2$ be homeomorphisms. Then
\[
 f(\s) \subset \bar{B} \left( h(\s), \| f-h \|_{\infty, \p\s} \right) .
\]
\end{lemma}

The following lemma will be useful to prove that our approximating piecewise affine function is a in fact homeomorphism.

\begin{lemma}\label{le:fs1fs2}
Let $\s_1, \s_2 \subset \R^2$ be two closed triangles such that $\s_1 \cap \s_2$ is empty, or a side of both $\s_1$ and $\s_2$, or a vertex of both $\s_1$ and $\s_2$.
Let $f : \s_1 \cup \s_2 \to \R^2$ be a function such that $f|_{\s_1}$ and $f|_{\s_2}$ are homeomorphisms, and $f$ is injective in $\p\s_1 \cup \p\s_2$. Suppose that for all vertices $w_1$ of the triangle $\s_1$ such that $w_1 \notin \s_2$ one has $f(w_1) \notin f(\s_2)$, and that for all vertices $w_2$ of the triangle $\s_2$ such that $w_2 \notin \s_1$ one has $f(w_2) \notin f(\s_1)$.
Then
\[
 f(\inte{\s}_1) \cap f(\s_2) = \varnothing .
\]
\end{lemma}
\begin{proof}
Suppose $f(\inte{\s}_1) \cap f(\s_2) \neq \varnothing$; we shall reach a contradiction.
So there exists $x\in \inte{\s}_1$ such that $f(x) \in f(\s_2)$.
In any of the three cases considered in the statement ($\s_1 \cap \s_2$ is empty, or a side of both $\s_1$ and $\s_2$, or a vertex of both $\s_1$ and $\s_2$), there exists a vertex $w$ of
$\s_1$ such that $w \notin \s_2$; by assumption, $f(w) \notin f(\s_2)$. Thanks to the geometry of the triangle $\s_1$, one has $[x,w) \subset \inte{\s}_1$. Define
\[
 t_0 := \sup \left\{ t \in [0,1] : f(x + t(w-x)) \in f(\s_2) \right\}.
\]
The number $t_0$ is well defined since $f(x) \in f(\s_2)$.
Moreover, $f(x + t_0 (w-x)) \in f(\s_2)$ because $f(\s_2)$ is closed.
As $f(w) \notin f(\s_2)$ then $t_0 < 1$ and, hence, $x + t_0 (w-x) \in \inte{\s}_1$.
Furthermore, by definition of supremum,
\[
 f(x + t(w-x)) \notin f(\s_2)  \qquad \mbox{for all} \quad t\in (t_0,1] .
\]
This proves that $f(x + t_0 (w-x)) \in \p f(\s_2)$. Since $f|_{\s_2}$ is a homeomorphism, by the invariance of domain Theorem, $\p f(\s_2) = f(\p \s_2)$. We have therefore proved $f(\inte{\s}_1) \cap f(\p \s_2) \neq \varnothing$.
Since $f|_{\s_1}$ is a homeomorphism, we have in fact
\[
 f(\inte{\s}_1) \cap f(\p \s_2 \setminus \p\s_1) \neq \varnothing.
\]
So there exists $y \in \p \s_2 \setminus \p\s_1$ such that $f(y) \in f(\s_1)$.
In any of the three cases considered in the statement ($\s_1 \cap \s_2$ is empty, or a side of both $\s_1$ and $\s_2$, or a vertex of both $\s_1$ and $\s_2$), we can find a vertex $v$ of the triangle $\s_2$ such that $v\notin \s_1$ and $[y,v] \subset \p \s_2 \setminus \p\s_1$. By assumption, $f(v) \notin f(\s_1)$.
Define
\[
 s_0 := \sup \left\{ s \in [0,1] : f(y + s(v-y)) \in f(\s_1) \right\}.
\]
The number $s_0$ is well defined since $f(y) \in f(\s_1)$.
Moreover, $f(y + s_0 (v-y)) \in f(\s_1)$ because $f(\s_1)$ is closed.
As $f(v) \notin f(\s_1)$ then $s_0 < 1$.
Furthermore, by definition of supremum,
\[
 f(y + s(v-y)) \notin f(\s_1)  \qquad \mbox{for all} \quad s\in (s_0,1] .
\]
This proves that $f(y + s_0 (v-y)) \in \p f(\s_1)$. Since $f|_{\s_1}$ is a homeomorphism, by the invariance of domain Theorem, $\p f(\s_1) = f(\p \s_1)$. We have therefore proved $f(\p\s_1) \cap f(\p \s_2 \setminus \p\s_1) \neq \varnothing$, and this is a contradiction with the fact that $f$ is injective in $\p\s_1 \cup \p\s_2$.
\end{proof}

The following is the main result of this section. It constructs a
piecewise linear homeomorphism that approximates a given
homeomorphism in the supremum norm. Its proof is based on the one by
Moise \cite{Moise} (Theorem 3 of Chapter 6) but here we also
estimate the side lengths and angles of the triangles of the
constructed triangulation.

\begin{theorem}\label{th:constructiontriangulation}
Let $\O \subset \R^2$ be a closed polygon. Let $\a, \tilde{\a} \in (0,1]$ and $H, \tilde{H}>0$. Then there exists
$\e_0>0$
such that for all $0< \e \leq \e_0$ and
for every homeomorphism $h \in C^{\a} (\O , \R^2)$ with $|h|_{\a} \leq H$ satisfying $h^{-1} \in C^{\tilde{\a}} (h(\O) , \R^2)$ and $|h^{-1}|_{\tilde{\a}} \leq \tilde{H}$,
there exist a triangulation $K$ of $\O$, and a homeomorphism $f: \O \to \R^2$ piecewise affine over $K$ such that $\| f-h \|_{\infty} \leq \e$,
\begin{align}
 \label{eq:inequality1} & \sin \f \geq A_0 \e^{a_0} \quad \mbox{for all angles } \ \f \ \mbox{ of all triangles in } K^2, \\
 \label{eq:inequality23} & A_1 \e^{a_1} \leq |e| \leq A_2 \e^{a_2} \quad \mbox{for all } \ e \in K^1 ,
\end{align}
where
\begin{equation}\label{eq:A}
\begin{split}
 & A_0 := C_0 \cdot B_1^{c_{01}} \cdot d^{\frac{b_1 c_{01}}{\tilde{\a}} + c_{02}} \cdot (\sin \t)^{\frac{b_1 c_{01}}{\tilde{\a}} + c_{03}} \cdot
 \tilde{H}^{-\frac{b_1 c_{01}}{\tilde{\a}}} \cdot (3H)^{-\frac{b_1 c_{01}}{\a \tilde{\a}} - \frac{c_{02}}{\a} + \frac{c_{04}}{\a}} , \\
 & A_1 := C_1 \cdot B_1^{c_{11}} \cdot d^{\frac{b_1 c_{11}}{\tilde{\a}} + c_{12}} \cdot (\sin \t)^{\frac{b_1 c_{11}}{\tilde{\a}} + c_{13}} \cdot \tilde{H}^{- \frac{b_1 c_{11}}{\a \tilde{\a}}} \cdot (3H)^{-\frac{b_1 c_{11}}{\tilde{\a}}-\frac{c_{12}}{\a} + \frac{c_{14}}{\a}} , \\
 & A_2 := (3H)^{-1/\a}, \\
 & a_0 := \frac{b_1 c_{01}}{\a \tilde{\a}} + \frac{c_{02}}{\a} - \frac{c_{04}}{\a} , \qquad
   a_1 := \frac{b_1 c_{11}}{\a \tilde{\a}} + \frac{c_{12}}{\a} - \frac{c_{14}}{\a} , \qquad
  a_2 := \frac{1}{\a},
\end{split}
\end{equation}
\eqref{eq:B}, \eqref{eq:td} and \eqref{eq:C0C1}.
\end{theorem}
\begin{proof}
Let $0 < \e \leq \e_0$, where $\e_0$ is to be decided later.
If $\e_0$ is small, by Proposition \ref{prop:triangulation}, there exist \eqref{eq:td}
and a triangulation $L$ of $\O$ such that
\begin{equation}\label{eq:diams}
 d \left( \frac{\e}{3H} \right)^{1/\a} \leq |e| \leq \left( \frac{\e}{3H} \right)^{1/\a}, \qquad e \in L^1,
\end{equation}
and \eqref{eq:singeq}.
By \eqref{eq:diams} and the geometry of the triangle,
\begin{equation}\label{eq:diamhs}
 \diam \s \leq \left( \frac{\e}{3H} \right)^{1/\a} , \qquad \diam h(\s) \leq \frac{\e}{3} , \qquad \s \in L^2 .
\end{equation}
By \eqref{eq:diams}, \eqref{eq:singeq} and the geometry of the complex $L$,
\begin{equation}\label{eq:lsint}
 \dist (v,\s) \geq d \sin \t \left( \frac{\e}{3H} \right)^{1/\a} , \qquad v\in \tilde{L}^0, \quad \s \in L^2, \quad v \notin \s.
\end{equation}
By \eqref{eq:lsint}, for all $\s \in L^2$,
\begin{equation}\label{eq:disthgeq}
 \dist (h(\s), h(\tilde{L}^0 \setminus \s)) \geq \left( \frac{\dist (\s, \tilde{L}^0 \setminus \s)}{\tilde{H}} \right)^{1/\tilde{\a}} \geq \left( \frac{d \sin \t}{\tilde{H}} \right)^{1/\tilde{\a}} \left( \frac{\e}{3H} \right)^{\frac{1}{\a\tilde{\a}}} .
\end{equation}
Call
\[
 \d:= \left( \frac{d \sin \t}{\tilde{H}} \right)^{1/\tilde{\a}} \left( \frac{\e}{3H} \right)^{\frac{1}{\a\tilde{\a}}} .
\]
If $\a \tilde{\a} < 1$ and $\e_0$ is small then
\begin{equation}\label{eq:dleqe3}
 \d \leq \frac{\e}{3} ;
\end{equation}
otherwise, if $\a = \tilde{\a} = 1$ then $H \tilde{H} \geq 1$ and \eqref{eq:dleqe3} as well.

By \eqref{eq:dleqe3}, \eqref{eq:singeq} and Theorem \ref{th:skeleton}, if $\e_0$ is small,
there exist a homeomorphism $g: \bigcup L^1 \to \R^2$
and a $1$-dimensional complex $M$ such that
\begin{equation}\label{eq:g-h}
 \bigcup M = \bigcup L^1 , \qquad L^0 \subset M^0 , \qquad \| g-h \|_{\infty, \bigcup M} < \d,
\end{equation}
$g$ is piecewise affine over $M$, coincides with $h$ in $\tilde{L}^0$, and
\begin{equation}\label{eq:|e|M}
 B_1 \d^{b_1} \leq |e| \leq B_2 \d^{b_2} , \qquad e \in M^1 ,
\end{equation}
where \eqref{eq:B}.

Let $\s \in L^2$. Note that $\tilde{M}^0 \cap \s \subset \p \s$. Thanks to Theorem \ref{th:Gupta} applied to
the triangle $\s$,
the complex
\[
 \{ e \in M^1 : e \subset \p \s \} \cup \left\{ \{ p \} \in M^0 : p \in \s \right\}
\]
and the homeomorphism $g|_{\p \s}$,
there exist a triangulation $N_{\s}$ of $\s$ and a homeomorphism $f_{\s} : \s \to \R^2$
piecewise affine over $N_{\s}$ that coincides with $g$ in $\p \s$ such that
\[
  \tilde{N}_{\s}^0 \cap \p \s = \tilde{M}^0 \cap \s,
\]
for all angles $\f$ of all triangles in $N_{\s}^2$ we have
\begin{equation}\label{eq:sintgeq}
 \sin \f \geq C_0 \left( B_1 \d^{b_1} \right)^{c_{01}} \left[ d \left( \frac{\e}{3H} \right)^{1/\a} \right]^{c_{02}} \left[ \left( \frac{\e}{3H} \right)^{1/\a} \right]^{-c_{04}} (\sin \t)^{c_{03}} ,
\end{equation}
and for all $e \in N_{\s}^1$ we have
\begin{equation}\label{eq:|e|N}
 C_1 \left( B_1 \d^{b_1} \right)^{c_{11}} \left[ d \left( \frac{\e}{3H} \right)^{1/\a} \right]^{c_{12}} \left[ \left( \frac{\e}{3H} \right)^{1/\a} \right]^{-c_{14}} (\sin \t)^{c_{13}} \leq |e| \leq \left( \frac{\e}{3H} \right)^{1/\a} ,
\end{equation}
where \eqref{eq:C0C1}.
In inequalities \eqref{eq:sintgeq} and \eqref{eq:|e|N} we have also applied \eqref{eq:diams}, \eqref{eq:singeq} and \eqref{eq:|e|M}.

Define $f: \O \to \R^2$ as the only function that coincides with $f_{\s}$ for each $\s \in L^2$.
Let $K$ be the $2$-dimensional complex defined by the condition
\[
 K^i = \bigcup_{\s \in L^2} N_{\s}^i , \qquad i = 0,1,2 .
\]
Then $K$ is a triangulation of $\O$, and $f$ is piecewise affine over $K$.
By construction and \eqref{eq:g-h}, $f$ coincides with $g$ in $\bigcup M = \bigcup L^1 = \bigcup K^1$, and $\| f-h \|_{\infty,\bigcup M} < \d$. Therefore, by Lemma \ref{le:fsubsetBh},
\begin{equation}\label{eq:fssubsetB}
 f(\s) \subset B (h(\s),\d), \qquad \s \in L^2 .
\end{equation}
Now, for each $x \in \O$ there exists $\s \in L^2$ such that $x \in \s$ and, hence,
\[
f(x) , h(x) \in B( h(\s) , \d) ;
\]
consequently, by \eqref{eq:diamhs} and \eqref{eq:dleqe3},
\[
 | f(x) - h(x) | < \diam B( h(\s) , \d) \leq \frac{\e}{3} + 2\d \leq \e .
\]
This proves $\| f-h \|_{\infty} < \e$.

For each $\D \in K^2$ there exists $\s \in L^2$ such that $\D \in N_{\s}^2$, and, hence, \eqref{eq:sintgeq} holds for all angles $\f$ of $\D$.
Similarly, for each $e \in K^1$ there exists $\s \in L^2$ such that $e \in N_{\s}^1$, and, hence,
\eqref{eq:|e|N} holds.

We now prove that $f$ is a homeomorphism. Since $f$ is continuous and $\O$ is compact, it suffices to show that $f$ is injective.
We have already seen that $f$ is injective in $\bigcup M$ (since so is $g$), and in $\s$ for each $\s \in L^2$ (since so is $f_{\s}$).
Now take $\s_1 \neq \s_2 \in L^2$.
As $L$ is a complex, then $\s_1 \cap \s_2$ is empty, or a side of both $\s_1$ and $\s_2$, or a vertex of both $\s_1$ and $\s_2$.
Take a vertex $w_1$ of $\s_1$ such that $w_1 \notin \s_2$. By \eqref{eq:disthgeq} and \eqref{eq:fssubsetB}, we have $f(w_1) = g(w_1) = h(w_1) \notin f(\s_2)$. Analogously, all vertices $w_2$ of $\s_2$ such that $w_2 \notin \s_1$ satisfy $f(w_2) \notin f(\s_1)$. By Lemma \ref{le:fs1fs2}, $f(\inte{\s}_1) \cap f(\s_2) = \varnothing$,
thus finishing the proof that $f$ is a homeomorphism.
\end{proof}

Note that, when $\a = \tilde{\a} = 1$, the coefficient $a_0$ of \eqref{eq:A} equals $0$ (by \eqref{eq:B} and \eqref{eq:C0C1}) and hence, because of \eqref{eq:inequality1}, the triangulation $K$ constructed in Theorem \ref{th:constructiontriangulation} is regular in the sense of Ciarlet \cite{Ciarlet}.

\section{Estimates in the H\"older norm}\label{se:estimates}

Following the notation of Theorem \ref{th:introduction}, up to now (Theorem \ref{th:constructiontriangulation}) we have constructed a piecewise affine homeomorphism $f$ that approximates $h$ in the supremum norm. In this section we will see how the approximation in the supremum norm and a control on the minimum and maximum side lengths and on the minimum angle of the triangulation (as done in Theorem \ref{th:constructiontriangulation}) will provide us with an approximation in the H\"older norm.

This section consists of three subsections.
In Subsection \ref{subse:trigonometricineq} we prove an elementary trigonometric inequality that will be useful for Subsection \ref{subse:apriori}.
In Subsection \ref{subse:apriori} we show a priori bounds in the H\"older norm of a piecewise affine function $u$ in terms of the  minimum and maximum side lengths and of the minimum angle of the triangulation over which $u$ is piecewise affine.
Subsection \ref{subse:approximation} uses all the results of the paper to prove Theorem \ref{th:introduction}.

\subsection{A trigonometric inequality}\label{subse:trigonometricineq}

In this subsection we show an elementary trigonometric inequality that will be used in Subsection \ref{subse:apriori}.
The following lemma plays a similar role to the one that Theorem 3.1.3 of Ciarlet \cite{Ciarlet} does in the context of proving approximation in the Sobolev norm for finite elements.
We use the standard notation that all elements of $\R^2$ are regarded as column vectors; in particular, given $a, b \in \R^2$, then $(a,b)$ is the $2 \times 2$ matrix whose columns are $a$ and $b$.
The norm $\|\cdot\|$ of a matrix is defined as the operator norm with respect to the Euclidean norm in $\R^2$.

\begin{lemma}\label{le:inradius}
Let $p_1, p_2, p_3 \in \R^2$ be three affinely independent points.
Then
\[
 \| (p_1 - p_3, p_2 - p_3)^{-1} \| \leq \frac{1}{\sin \an{p_1 p_3 p_2}} \left( \frac{1}{|p_1 - p_3|} + \frac{1}{|p_2 - p_3|} \right) .
\]
\end{lemma}
\begin{proof}
Elementary matrix computations show that
\begin{align*}
 \left\| (p_1 - p_3, p_2 - p_3)^{-1} \right\| & \leq \left\| (p_1 - p_3, p_2 - p_3)^{-1} \right\|_2 = \frac{\left\| (p_1 - p_3, p_2 - p_3) \right\|_2}{|\det (p_1 - p_3, p_2 - p_3)|}  \\
 & \leq  \frac{|p_1 - p_3| +  |p_2 - p_3|}{|\det (p_1 - p_3, p_2 - p_3)|} ,
\end{align*}
where by definition $\|\cdot\|_2$ is the Frobenius norm of a matrix, that is, the Euclidean norm of a matrix regarded as a vector in $\R^4$.
By elementary plane geometry,
\[
 |\det (p_1 - p_3, p_2 - p_3)| = 2 \area \left( p_1 p_2 p_3 \right) = |p_1 - p_3| |p_2 - p_3| \sin \an{p_1 p_3 p_2} ,
\]
which concludes the proof.
\end{proof}

\subsection{A priori bounds in the H\"older norm of piecewise affine functions}\label{subse:apriori}

In this subsection we prove two a priori bounds in the H\"older norm of piecewise affine functions, in terms of the minimum and maximum side length and on the minimum angle of the triangulation.

First we recall a standard geometric property of Lipschitz domains (see, e.g., Exercise 1.9 of Ciarlet \cite{CiarletElasticity}).

\begin{lemma}\label{le:cO}
Let $\O \subset \R^2$ be the closure of a non-empty, open, bounded, connected set with Lipschitz boundary. Then there exists a constant
\begin{equation}\label{eq:cO}
 c(\O)\geq 1 \quad \mbox{depending only on } \ \O 
\end{equation}
such that for each $x,y \in \O$ there exist
\begin{equation}\label{eq:mqq}
 m\geq 1 \quad \mbox{and} \quad q_0, \ldots, q_m \in \R^2
\end{equation}
satisfying 
\begin{equation}\label{eq:conditionscO}
  q_0 = x, \qquad q_m = y, \qquad \bigcup_{i=0}^{m-1} [q_i, q_{i+1}] \subset \O , \qquad
  \sum_{i=0}^{m-1} |q_i - q_{i+1}| \leq c(\O) |x-y| .
\end{equation}
\end{lemma}

We will also need the following standard and straightforward interpolation inequality between H\"older spaces.

\begin{lemma}\label{le:interpolation}
Let $\O \subset \R^2$.
Consider $0 < \b \leq \a \leq 1$ and $u \in C^{\a} (\O, \R^2)$.
Then
\[
 |u|_{\b} \leq 2^{1-\frac{\b}{\a}} \|u\|_{\infty}^{1-\frac{\b}{\a}} |u|_{\a}^{\frac{\b}{\a}} .
\]
\end{lemma}

We will use Lemma \ref{le:interpolation} in the following way.
Following the notation of Theorem \ref{th:constructiontriangulation}, we have proved that $\|f-h\|_{\infty}$ is small.
Therefore, by Lemma \ref{le:interpolation}, to prove that $\|f-h\|_{\b}$ is small, we only have to show a priori bounds on $|f-h|_{\a}$. These bounds will be calculated in Propositions \ref{prop:piecewise} and \ref{prop:errorinterpolant}. In fact, we will bound $|f-h|_{\a}$ in terms of a negative power of $\e$ (recall from Theorem \ref{th:constructiontriangulation} that $\|f-h\|_{\infty} \leq \e$), but we will see in Subsection \ref{subse:approximation} that this suffices to obtain Theorem \ref{th:introduction}.

\begin{proposition}\label{prop:piecewise}
Let $\O \subset \R^2$ be a closed polygon. Consider real numbers
\begin{equation}\label{eq:eacccaaa}
 0 < \a \leq 1, \quad
 0 < A_0, A_1, A_2, \quad
 0 \leq a_0, \quad
 a_1 \geq a_2 .
\end{equation}
Then there exist $c_3 > 0$ such that for every $\e >0$ and
every $u \in C (\O , \R^2)$ piecewise affine over any triangulation $K$
of $\O$ satisfying \eqref{eq:inequality1} and \eqref{eq:inequality23}, we have
\begin{equation}\label{eq:ualequinfty}
 | u |_{\a} \leq c_3 \|u\|_{\infty} \e^{a_3} ,
\end{equation}
where
\begin{equation}\label{eq:a4}
 a_3 := -\a(a_0 + a_1) .
\end{equation}
\end{proposition}
\begin{proof}
Take $\s \in K^2$ and $x,y \in \s$. 
Let $p_1, p_2, p_3$ be the three vertices of the triangle $\s$.
It is easy to see that, when we define
\begin{equation}\label{eq:lmu}
 \binom{\l_x}{\mu_x} = (p_1 - p_3, p_2 - p_3)^{-1} (x-p_3) , \qquad
 \binom{\l_y}{\mu_y} = (p_1 - p_3, p_2 - p_3)^{-1} (y-p_3) ,
\end{equation}
then we have
\begin{equation*}
\begin{split}
  x & = \l_x p_1 + \mu_x p_2 + (1- \l_x - \mu_x) p_3 , \\
  y & = \l_y p_1 + \mu_y p_2 + (1- \l_y - \mu_y) p_3 ;
\end{split}
\end{equation*}
consequently, since $u|_{\s}$ is affine,
\begin{equation}\label{eq:formulafxfy}
\begin{split}
  u(x) & = \l_x u(p_1) + \mu_x u(p_2) + (1- \l_x - \mu_x) u(p_3) , \\
  u(y) & = \l_y u(p_1) + \mu_y u(p_2) + (1- \l_y - \mu_y) u(p_3) .
\end{split}
\end{equation}
By \eqref{eq:lmu},
\begin{equation}\label{eq:l-mu}
 \binom{\l_x - \l_y}{\mu_x - \mu_y} = (p_1 - p_3, p_2 - p_3)^{-1} (x-y) .
\end{equation}
From \eqref{eq:formulafxfy} and \eqref{eq:l-mu} we obtain
\begin{equation}\label{eq:fx-fy}
 u(x) - u(y) = \left( u(p_1) - u(p_3) , u(p_2) - u(p_3) \right) (p_1 - p_3, p_2 - p_3)^{-1} (x-y) .
\end{equation}
Now we use Lemma \ref{le:inradius} and inequalities \eqref{eq:inequality1}, \eqref{eq:inequality23} to obtain
\begin{equation}\label{eq:()-1leqe}
 \left\|(p_1 - p_3, p_2 - p_3)^{-1}\right\| \leq
 2 A_0^{-1} A_1^{-1} \e^{- a_0 - a_1}.
\end{equation}
Equations \eqref{eq:fx-fy} and \eqref{eq:()-1leqe} show that
\begin{equation}\label{eq:basic1}
  |u(x)-u(y)| \leq 4 \sqrt{2} \|u\|_{\infty} A_0^{-1} A_1^{-1} \e^{- a_0 - a_1} |x-y| .
\end{equation}

Now let $x,y \in \O$ be arbitrary.
By Lemma \ref{le:cO}, there exist \eqref{eq:cO} and \eqref{eq:mqq} such that \eqref{eq:conditionscO}.
For each $i \in \{0, \ldots, m-1\}$ let $m_i \geq 1$ and $q_{i,0}, \ldots, q_{i,m_i} \in \R^2$ be such that
\begin{align*}
 & q_{i,0} = q_i, \qquad q_{i,m_i} = q_{i+1}, \qquad 
 [q_{i,j}, q_{i,j+1}] \subset \s_j \quad \mbox{for some } \ \s_j \subset K^2 , \\
 & q_{i,0} < \cdots < q_{i,m_i} \quad \mbox{in the order of } \ [q_i, q_{i+1}] .
\end{align*}
Then, by \eqref{eq:basic1},
\[ |u(x) - u(y)|  \leq \sum_{i=0}^{m-1} \sum_{j=0}^{m_i-1} |u(q_{i,j}) - u(q_{i,j+1})| \leq 4 \sqrt{2} \|u\|_{\infty} A_0^{-1} A_1^{-1} \e^{- a_0 - a_1} \sum_{i=0}^{m-1} \sum_{j=0}^{m_i-1} |q_{i,j} - q_{i,j+1}|.
\]
But now, by \eqref{eq:conditionscO},
\[
 \sum_{i=0}^{m-1} \sum_{j=0}^{m_i-1} |q_{i,j} - q_{i,j+1}| = \sum_{i=0}^{m-1} |q_i - q_{i+1}| \leq c(\O) |x-y| .
\]
This proves
$|u|_1 \leq 4 \sqrt{2} \|u\|_{\infty} A_0^{-1} A_1^{-1} c(\O) \e^{- a_0 - a_1}$,
and, hence, by Lemma \ref{le:interpolation},
\[
 |u|_{\a} \leq 2^{1-\a} \|u\|_{\infty}^{1-\a} |u|_1^{\a} \leq 
  2^{1+ \frac{3}{2}\a} A_0^{-\a} A_1^{-\a} c(\O)^{\a} \|u\|_{\infty} \e^{-\a(a_0 + a_1)} .
\]
This concludes the proof.
\end{proof}

\begin{proposition}\label{prop:errorinterpolant}
Let $\O \subset \R^2$ be a closed polygon. Let \eqref{eq:eacccaaa}.
Then there exist $c_4 > 0$ such that for every $\e >0$,
every triangulation $K$ of $\O$ satisfying \eqref{eq:inequality1} and \eqref{eq:inequality23},
and every $h \in C^{\a} (\O , \R^2)$, we have
\begin{equation}\label{eq:h-Pih}
 | \Pi h |_{\a} \leq c_4 \| h \|_{\a} \e^{a_4} ,
\end{equation}
where
\begin{equation}\label{eq:a5}
 a_4 := -\a (a_0 + a_1 - \a a_2)
\end{equation}
and $\Pi h$ is the piecewise affine function over $K$ that coincides with $h$ in $\tilde{K}^0$.
\end{proposition}
\begin{proof}
Take $\s \in K^2$ and $x,y \in \s$. 
Let $p_1, p_2, p_3$ be the three vertices of the triangle $\s$.
By \eqref{eq:fx-fy},
\begin{equation}\label{eq:Pihx-Pihy}
 \Pi h(x) - \Pi h(y) = \left( h(p_1) - h(p_3) , h(p_2) - h(p_3) \right) (p_1 - p_3, p_2 - p_3)^{-1} (x-y) .
\end{equation}
Now, thanks to \eqref{eq:inequality23},
\begin{equation}\label{eq:case1h}
 \left\| ( h(p_1) - h(p_3) , h(p_2) - h(p_3) ) \right\| \leq |h|_{\a} \left( |p_1 - p_3|^{\a} + |p_2 - p_3|^{\a} \right)
 \leq 2 |h|_{\a} A_2^{\a} \e^{\a a_2} ;
\end{equation}
on the other hand, by Lemma \ref{le:inradius}, \eqref{eq:inequality1} and \eqref{eq:inequality23},
\begin{equation}\label{eq:case1p}
 \| (p_1 - p_3, p_2 - p_3)^{-1} \| \leq \frac{1}{\sin \an{p_1 p_3 p_2}} \left( |p_1 - p_3|^{-1} + |p_2 - p_3|^{-1} \right) \leq 2 A_0^{-1} A_1^{-1} \e^{-a_0-a_1}.
\end{equation}
In total, \eqref{eq:Pihx-Pihy}, \eqref{eq:case1h} and \eqref{eq:case1p} show that
\begin{equation}\label{eq:basic1a}
 |\Pi h(x) - \Pi h(y)| \leq 4 |h|_{\a} A_0^{-1} A_1^{-1} A_2^{\a}  \e^{-a_0-a_1+\a a_2}  |x-y| .
\end{equation}

Now, using Lemma \ref{le:cO} and \eqref{eq:basic1a}, and arguing as in the proof of Proposition \ref{prop:piecewise}, one can show easily that
\begin{equation}\label{Pih1}
 |\Pi h|_1 \leq 4 |h|_{\a} A_0^{-1} A_1^{-1} A_2^{\a} c(\O) \e^{-a_0-a_1+\a a_2} ,
\end{equation}
where $c(\O)$ is the constant of Lemma \ref{le:cO}.
Now we use \eqref{Pih1}, Lemma \ref{le:interpolation} and the obvious inequalities
$\|\Pi h\|_{\infty} \leq \|h\|_{\infty} \leq \|h\|_{\a}$ and $|h|_{\a} \leq \|h\|_{\a}$
to conclude that 
\[
  |\Pi h|_{\a} \leq 2^{1+\a} A_0^{-\a} A_1^{-\a} A_2^{\a^2} c(\O)^{\a} \|h\|_{\a} \e^{\a(-a_0-a_1+\a a_2)} ,
\]
thus finishing the proof.
\end{proof}

\subsection{Approximation in the H\"older norm}\label{subse:approximation}
Up to now (Theorem \ref{th:constructiontriangulation}) we have constructed a piecewise affine homeomorphism $f$ that approximates $h$ in the supremum norm.
The a priori bounds found in Subsection \ref{subse:apriori} will provide us with an approximation in the H\"older norm, thus proving Theorem \ref{th:introduction}.

Recall that the exponents $a_0, a_1, a_2, a_3, a_4$ defined in \eqref{eq:A}, \eqref{eq:a4} and \eqref{eq:a5} satisfy
\begin{equation}\label{eq:aineq}
 0 \leq a_0, \qquad 1 \leq a_2 \leq a_1 , \qquad a_4 \leq 1 + a_3 \leq 0 , 
\end{equation}
and that these inequalities are equalities if $\a = \tilde{\a} = 1$.

Theorem \ref{th:constructiontriangulation}, Propositions \ref{prop:piecewise} and \ref{prop:errorinterpolant} and Lemma \ref{le:interpolation} provide a proof of the main theorem of this paper, stated below.

\begin{theorem}\label{th:main}
Let $\O \subset \R^2$ be a closed polygon. Let $0 < \a, \tilde{\a} \leq 1$ and let $h\in C^{\a} (\O,\R^2)$ be a homeomorphism with $h^{-1} \in C^{\tilde{\a}} (h(\O),\R^2)$.
Then there exist constants
$\e_0 , D >0$
depending only on 
\begin{equation}\label{eq:depending}
 \O, \ \|h\|_{\a}, \ |h^{-1}|_{\tilde{\a}} , \ \a, \ \tilde{\a}
\end{equation}
such that for every
\begin{equation}\label{eq:eandbmain}
 0 < \e \leq \e_0 \quad \mbox{and} \quad 0 < \b < \frac{\a}{ 1- a_4 }
\end{equation}
(where $a_4$ is defined through \eqref{eq:a5}, \eqref{eq:A}, \eqref{eq:C0C1} and \eqref{eq:B})
there exists a piecewise affine homeomorphism $f \in C (\O,\R^2)$ such that
\[
 \| f-h \|_{\infty} \leq \e \quad \mbox{and} \quad | f-h |_{\b} \leq D \e^{1 - \frac{\b}{\a} (1- a_4)} .
\]
\end{theorem}
\begin{proof}
Take \eqref{eq:eandbmain}, where $\e_0$ is to be chosen later.
If $\e_0$ is small, by Theorem \ref{th:constructiontriangulation},
there exist a triangulation $K$ of $\O$, and a homeomorphism $f: \O \to \R^2$ such that
\eqref{eq:inequality1}, \eqref{eq:inequality23},
$f$ is piecewise affine over $K$, and
\begin{equation}\label{eq:f-ge}
 \| f-h \|_{\infty} \leq \e.
\end{equation}
By Lemma \ref{le:interpolation}, we have 
\begin{equation}\label{eq:f-hb}
 |f-h|_{\b} \leq 2^{1-\frac{\b}{\a}} \|f-h\|_{\infty}^{1-\frac{\b}{\a}} |f-h|_{\a}^{\frac{\b}{\a}} .
\end{equation}
Now,
\begin{equation}\label{eq:hPihf}
 |f-h|_{\a} \leq |f- \Pi h|_{\a} + |\Pi h - h|_{\a} ,
\end{equation}
where $\Pi h$ is the piecewise affine function over $K$ that coincides with $h$ in $\tilde{K}^0$.
If $\e_0$ is small, by Proposition \ref{prop:piecewise} and Theorem \ref{th:constructiontriangulation}, there exists a constant $c_3 >0$ depending only on \eqref{eq:depending} such that
\[
 | f - \Pi h|_{\a} \leq c_3 \|f - \Pi h\|_{\infty} \e^{a_3} .
\]
We now note that $\|f - \Pi h\|_{\infty}\leq \e$; indeed, this is immediate, since $f - \Pi h$ is a piecewise affine function over $K$ such that (thanks to \eqref{eq:f-ge}) $|f(x) - \Pi h(x)| \leq \e$ for all $x \in \tilde{K}^0$. Therefore,
\begin{equation}\label{fPihe}
 | f - \Pi h|_{\a} \leq c_3 \e^{1+a_3} .
\end{equation}
By Proposition \ref{prop:errorinterpolant} and Theorem \ref{th:constructiontriangulation} there exists a constant $c_4\geq 1$ depending only on \eqref{eq:depending} such that inequality \eqref{eq:h-Pih} holds, and hence,
\begin{equation}\label{eq:pih-h}
 |\Pi h - h|_{\a} \leq |\Pi h|_{\a} + |h|_{\a} \leq c_4 \|h\|_{\a} \e^{a_4} + \|h\|_{\a} \leq 2 c_4 \|h\|_{\a} \e^{a_4} ,
\end{equation}
since $\e_0$ is small and (by \eqref{eq:aineq}) $a_4 \leq 0$.
In total, inequalities \eqref{eq:f-ge}, \eqref{eq:f-hb}, \eqref{eq:hPihf}, \eqref{fPihe}, \eqref{eq:pih-h} and \eqref{eq:aineq} demonstrate that
\[
 |f-h|_{\b} \leq 2^{1-\frac{\b}{\a}} \left( c_3 + 2 c_4 \|h\|_{\a} \right)^{\frac{\b}{\a}} \e^{1 - \frac{\b}{\a} (1-a_4)} .
\]
This concludes the proof.
\end{proof}

Of course, Theorem \ref{th:main} demonstrates Theorem \ref{th:introduction}.

We finish this paper with some comments about the optimality of Theorem \ref{th:main}.

Let $0 < \a, \tilde{\a} \leq 1$. Let $B(\a,\tilde{\a})$ be the set of all $0 < \b \leq 1$ with the following property:
For every closed polygon $\O$, the set of piecewise affine homeomorphisms from $\O$ to $\R^2$ is dense in the set
\[
 \left\{ h \in C^{\a} (\O,\R^2) : h \mbox{ is a homeomorphism with } h^{-1} \in C^{\tilde{\a}}(h(\O),\R^2) \right\}  
\]
in the $\|\cdot\|_{\b}$ norm.
Clearly, $B(\a,\tilde{\a})$ is an interval.
What Theorem \ref{th:main} asserts is that $B(\a,\tilde{\a})$ is non-empty, and
\[
 \sup B(\a,\tilde{\a}) \geq \frac{\a}{1-a_4} = \frac{\a^4 \tilde{\a}^3}{3 + 6 \a \tilde{\a} - 6\a^2 \tilde{\a}^2 - \a^3\tilde{\a}^3 - \a^4\tilde{\a}^3} 
\]
(the latter equality comes from \eqref{eq:a5}, \eqref{eq:A}, \eqref{eq:C0C1} and \eqref{eq:B}).
An optimal version of the result of Theorem \ref{th:main} would be to calculate $\sup B(\a,\tilde{\a})$ and to ascertain whether $\sup B(\a,\tilde{\a})$ belongs to $B(\a,\tilde{\a})$.
For example, it is easy to show that
\begin{equation}\label{eq:anotin}
 \a \notin B(\a,\tilde{\a}) , \qquad \a < 1, \quad \tilde{\a} \leq 1 .
\end{equation}
Indeed, take $0 <\a <1$, and let $\O$ be a closed polygon whose interior contains $0$.
Define $h: \O \to \R^2$ by
\[
 h(x,y) = (|x|^{\a} \sgn x , y) , \qquad (x,y) \in \O ,
\]
where $\sgn x$ stands for the sign of $x\in\R$.
Then $h \in C^{\a} (\O,\R^2)$ is a homeomorphism, and $h^{-1}$ is Lipschitz continuous.
However, $h$ cannot be approximated by piecewise affine homeomorphisms in the $\|\cdot\|_{\a}$ norm.
In fact, it is easy to see that $h$ cannot be approximated by Lipschitz continuous functions in the $\|\cdot\|_{\a}$ norm (see, e.g., Kichenassamy \cite{Kichenassamy}, if necessary).
Therefore, $\a \notin B(\a,1)$, and this implies \eqref{eq:anotin}.

Note that Theorem \ref{th:main} shows in particular that $\sup B(1,1) = 1$, but we do not know whether $1$ belongs to $B(1,1)$.
In any case, we believe that Theorem \ref{th:main} is not optimal.

\begin{small}
\subsection*{Acknowledgements}
We thank John Ball for helpful and stimulating discussions. J. C. B.
acknowledges financial support from MEC (Spain), grant
MTM-2004-07114, and Junta de Comunidades de Castilla-La Mancha
(Spain), grant PAI-05-027. C. M.-C. has been financially supported by a fellowship from `Secretar\'ia de Estado de Educaci\'on y Universidades' (Spain) and co-funded by the European Social Fund.

\end{small}
\end{document}